\documentclass[11pt,a4paper]{article}
\usepackage[utf8]{inputenc}
\usepackage[english]{babel}
\usepackage{amsthm}
\usepackage{float}
\usepackage{subfig}
\usepackage{xr}
\usepackage{bm}
\usepackage{mathtools}
\usepackage{indentfirst}
\usepackage[T1]{fontenc}
\usepackage{color}
\usepackage{amsfonts} 
\usepackage{amsmath}
\usepackage{pgfplots}
\usepackage{tikz}
\usepackage{dsfont}
\usepackage[hyperindex=true,colorlinks=true]{hyperref}

\usepackage{mathrsfs}

\usepackage{pstricks}
\usepackage{esint}
\tikzset{math3d/.style={x={(-0.353cm,-0.353cm)},z={(0cm,1cm)},y={(1cm,0cm)}}}
\usetikzlibrary{calc}
\usetikzlibrary{patterns}
\usetikzlibrary{arrows} 

\newtheorem{Theo}{Theorem}[section]

\newtheorem{Prop}{Proposition}[section]
\newtheorem{Lem}{Lemma}[section]

\newtheorem{Def}{Definition}[section]
\newtheorem{Rem}{Remark}[section]

\DeclareMathOperator{\divg}{div}

\newcommand{\w}{\widetilde}
\newcommand{\R}{\mathbb{R}}
\newcommand{\dx}{\ \textrm{d}x}
\newcommand{\dy}{\ \textrm{d}y}

\newcommand{\harpoon}{\ -\hspace{-0.2cm}\rightharpoonup \  }
\newcommand{\T}{\textrm{T}}

\newcommand{\e}{\varepsilon}
\newcommand{\N}{\mathbb N}
\newcommand{\Z}{\mathbb Z}
\newcommand{\n}{\nabla}
\renewcommand{\qed}{\hfill \ensuremath{\Box}}
\newcommand{\supp}{\operatorname{supp}}

\newcommand{\1}{\mathds{1}_{\Omega\setminus \Omega_n}}
\newcommand{\2}{\mathds{1}_{\Omega_n}}

\def\Sum{\displaystyle\sum}

\tikzset { xmin/.store in=\xmin, xmin/.default=-3, xmin=-3,
           xmax/.store in=\xmax, xmax/.default=3, xmax=3,
           ymin/.store in=\ymin, ymin/.default=-3, ymin=-3,
           ymax/.store in=\ymax, ymax/.default=3, ymax=3}

\tikzset {domaine/.style 2 args={domain=#1:#2}}

\makeatletter
\renewcommand{\thefigure}{\ifnum \c@section>\z@ \thesection.\fi
 \@arabic\c@figure}
\@addtoreset{figure}{section}
\makeatother
\makeatletter
\renewcommand{\theequation}{\ifnum \c@section>\z@ \thesection.\fi
 \@arabic\c@equation}
\@addtoreset{equation}{section}
\makeatother

\title{\LARGE  Homogenization of high-contrast and non symmetric conductivities for non periodic columnar structures}
\author{\begin{tabular}{cc}
    {M. CAMAR-EDDINE } & {L. PATER}
    \\*[-.0em]
    {\small  Institut de Recherche Math\'ematique de Rennes} & {\small Institut de Recherche Math\'ematique de Rennes}
    \\*[-.3em]
    {\small  INSA de Rennes} &	 {\small Universit\'e de Rennes 1}
    \\*[-.3em]
    {\small camar@insa-rennes.fr} & {\small laurent.pater@ens-cachan.org}
\end{tabular}}

\begin{document}
\maketitle


\label{chap2}

\markboth{Chapitre~\ref*{chap2} - Homogenization of high-contrast and non symmetric conductivities for non periodic columnar structures}{\textit{Section~\ref*{intro} - Introduction}}

\begin{abstract}
In this paper we determine, in dimension three, the effective conductivities of non periodic high-contrast two-phase cylindrical composites, placed in a constant magnetic field, without any assumption on the geometry of their cross sections. Our method, in the spirit of the \mbox{H-convergence} of Murat-Tartar, is based on a compactness result and the cylindrical nature of the microstructure. The homogenized laws we obtain extend those of the periodic fibre-reinforcing case of~\cite{BP1} to the case of periodic and non periodic composites with more general transversal geometries. 
\end{abstract}


\section{Introduction}\label{intro}

\markboth{Chapitre~\ref*{chap2} - Homogenization of high-contrast and non symmetric conductivities for non periodic columnar structures}{\textit{Section~\ref*{intro} - Introduction}}

\noindent
At the end of the~${19}$th century, it was discovered~\cite{Hall} that a constant magnetic field~${h}$ modifies the symmetric conductivity matrix~${\sigma}$ of a conductor into a non symmetric matrix~${\sigma(h)}$. This is know as the Hall effect. In the Maclaurin series of the perturbed resistivity~$(\sigma(h))^{-1}$ the zeroth-order term coincides with the resistivity~${\sigma^{-1}}$ in the absence of a magnetic field~\cite{LaLi}. In dimension two,~$h$ is a scalar and the first-order term is an antisymmetric matrix proportional to~${hJ}$; the coefficient of proportionality is called the Hall coefficient.  In dimension three,~$h \in \R^3$ and the first-order term, in the Maclaurin series of~${(\sigma(h))^{-1}}$, is of the form~${\mathscr E(R h)}$ where~${\mathscr E(\xi)j := \xi \times j}$ and~${R}$ is a~${3\times3}$ matrix called the Hall matrix~\cite{BM1}. 
In this work, we consider the idealized situation when the induced non symmetric part is proportional to the applied magnetic field:~${\sigma(h)=\alpha I_3+\beta{\mathscr E}(h)}$, where~${\alpha}$ and~${\beta}$ are two constant real numbers. For a given sequence of perturbed conductivities~${\sigma_n(h)}$, it is of great interest, in electrodynamics~\mbox{\cite{LaLi, Ali}}, to understand the influence of the magnetic field~${h}$ on the effective Hall coefficient or the effective Hall matrix through the homogenization of~${\sigma_n(h)}$. 

\medskip
Let us first review a few of the mathematical theory of  homogenization of elliptic partial differential equations of the form
\begin{equation} \label{eqdiffintro}
 \left\{\begin{array}{rcll}
- \divg \big( \sigma_n \n u_n \big) & = & f & \text{in } \Omega, \\*[.4em]
u_n & = & 0 & \text{on } \partial \Omega,
\end{array}\right.
\end{equation}
where~${\Omega}$ is a bounded open subset of~${\R^3}$,~${\sigma_n}$ is a sequence of matrix-valued functions in~${L^\infty(\Omega)^{3 \times 3}}$ and~${f}$ is an element of~${H^{-1}(\Omega)}$. This topic  has been intensively studied for the last four decades providing a wide literature~\cite{Spa,Mur1,Mur2,BLP}. In the context of conduction, when the conductivity matrices~${\sigma_n}$ are uniformly bounded, Spagnolo~\cite{Spa} with the~${G}$-convergence theory,  Murat and \mbox{Tartar}~\mbox{\cite{Mur1, Mur2}} with the~${H}$-convergence theory showed that the solution~${u_n\in H_0^1(\Omega)}$ of the conductivity problem~\eqref{eqdiffintro} strongly converges in~${L^2(\Omega)}$, up to a subsequence of~${n}$, to the solution of a limit conductivity problem of the same nature. The case of high-contrast conductivities is very different since non classical  phenomena, such as nonlocal terms, may appear in the limit problem as shown, for instance, in~\cite{FeKh,Khr,BelBou,BrTc,Bri1,KhMa}. This does not happen in dimension two if the sequence~${\sigma_n}$ is uniformly bounded from below. Briane~\cite{Brinonloc} and Casado-D\'ias \& Briane~\cite{BrCa1} proved that in that case the class of equations~(\ref{eqdiffintro}) is always compact in the sense that the limit equation of~\eqref{eqdiffintro} is always of the same type. In~\cite{BrCa1} they proved some extensions of the  well-known div-curl lemma of Murat-Tartar~\cite{Mur2} and deduce several compactness results under the assumption of  equicoerciveness coupled with the~${L^1}$-boundedness of the sequence of conductivities. 

\medskip
In this paper we are interested in the homogenization of a class of three-dimensional conductivity problems of the type
\begin{equation} \label{eqdiffnhintro} 
\left\{\begin{array}{r c l l}
- \divg \big( \sigma_n(h) \n u_n \big) & = & f & \text{in } \Omega, \\*[.4em]
u_n & = & 0 & \text{on } \partial \Omega,
\end{array}\right.
\end{equation}
where~${\sigma_n(h)}$ is an equi-coercive sequence of high-contrast two-phase conductivities perturbed by a constant magnetic field~${h\in\R^3}$ of the form~${\sigma_n(h):=(1-\mathds{1}_{{\Omega}_n})\sigma_{1}(h)+\mathds{1}_{{\Omega}_n}\sigma_{2,n}(h)}$ where~${\sigma_{2,n}(h)}$ is the unbounded perturbed conductivity of the highly conducting phase~${\Omega_n}$ and~${\sigma_{1}(h)}$ is the perturbed  conductivity of the phase surrounding~${\Omega_n}$. 

\medskip
In dimension two, for the case of low magnetic field, Bergman~\cite{Ber} was the first author who came up with a general formula for the effective Hall coefficient of a periodic composite material in terms of the local Hall coefficients and some local currents solving the conductivity equations in the absence of a magnetic field. We refer also to the works~\cite{MilHall, BerLiSt, BMM} for other two-dimensional composites, to~\mbox{\cite{BerSt2,BerSt1,StrelBer,Grab1,Grab2}} for composites with microstructure independent of one coordinate~(the so-called columnar composites) and to~\mbox{\cite{BerSt1a,BerSS}} for the case of strong magnetic field.

\medskip
Recently, in dimension two, M. Briane and the second author~\cite{BP1} obtained the effective perturbed conductivity~${\sigma_\ast(h)}$ of a sequence of isotropic high-contrast two-phase conductivities~${\sigma_n(h)}$ in the case of strong magnetic field, \emph{i.e.,} when the symmetric part and the antisymmetric part of the conductivity are of the same order. By extending a duality principle from~\cite{NHM} and using a suitable Dykhne transformation, which~(following Milton~\mbox{\cite{MilHall,MilBook}})  changes non symmetric matrices into symmetric ones, they proved that the symmetric part of the effective perturbed conductivity~${\sigma_\ast(h)}$ is given in terms of the effective conductivity in the absence of a magnetic field. They subsequently compared their two-dimensional results   to a three-dimensional periodic one and showed that the way a magnetic field perturbs the conductivity of a composite depends on the dimension. In order to compute the explicit perturbation formula in dimension three, they restricted themselves to a particular periodic fibre-reinforced structure, \emph{i.e.,} a structure completely described by any two-dimensional cross section transversal to the fibres~(first introduced by Fenchenko, Khruslov~\cite{FeKh} to derive a non local effect in homogenization). To our knowledge, only few results are known on the homogenization of both high-contrast and non symmetric conductivities in dimension three. 

\medskip
The aim of this paper is to determine the effective perturbed conductivity of~\eqref{eqdiffnhintro} for non periodic high-contrast two-phase cylindrical composites without any assumption on the geometry of the transversal microstructure. 

\medskip

We first investigate the periodic case, that is, when~${\sigma_n(h)(\cdot)= \Sigma_n(h)(\cdot /\e_n)}$ where~${\Sigma_n(h)(\cdot)}$ is a~\mbox{$Y$-periodic} matrix-valued function and~${\varepsilon_n \to 0}$ represents the size of the heterogeneities in the composite. 
In order to avoid non local effects in the limit problem, following Briane~\cite{BrARMA}, we assume the existence of a sequence of positive numbers~${c_n}$ such that~${\varepsilon_n^2 \, c_n}$ tends to zero, as~${n}$ goes to infinity, and satisfying the weighted Poincar\'e-Wirtinger inequality 
\begin{equation*} 
\forall \, V \in H^1(Y), \quad\int_{Y} |\Sigma_n(h)(y)| \left| V-\int_Y V\,\text{d}y\right|^2\,\text{d}y \leq c_n \int_{Y} \Sigma_n(y) \nabla V \cdot \n V\,\text{d}y.
\end{equation*} 
For a fixed~${n \in \mathbb N^*}$, using the theory of exact relations of Grabovsky, Milton, Sage~\mbox{\cite{Grab1,GMS}}~(thanks to the independence of the microstructure of the variable~${x_3}$), we obtain the H-limit~${\big( \sigma_n \big)_*}$ associated with the periodic homogenization~\cite{BLP} of the oscillating sequence~${\Sigma_n(\cdot /\e)}$ as~${\e \to 0}$. Then, we show that the sequence of constant conductivities~${\big( \sigma_n \big)_*}$ converges to some $\sigma_*(h)$ which, according to~\cite{BrARMA}, coincides with the homogenized conductivity associated with the limit problem of~\eqref{eqdiffnhintro}. The obtained effective conductivity~${\sigma_*(h)}$ is explicitly computed in terms of the homogenized conductivity~${\widetilde\sigma_\ast(h)}$ of the conduction problem posed in the~${(x_1,x_2)}$-plane transversal to the columnar composite~(see Proposition~\ref{propper}).
 
\medskip
 
Most of the arguments and tools used in the periodic case crucially lie on the periodic nature of the microstructure. 
 Therefore, a fundamentally different approach is necessary for the analysis of~\eqref{eqdiffnhintro} when~${\sigma_n(h)}$ is not periodic.

\medskip

In order to study the asymptotic behavior of the problem \eqref{eqdiffnhintro} in the non periodic case, using a method, in the spirit of the \mbox{H-convergence} of Murat-Tartar, we determine the limit, in an appropriate sense, of the current~${\sigma_n(h) \n u_n}$. The key ingredient of this approach is a fundamental compactness result~(see Lemma~\ref{divcurl}) based on a control of high conductivities in thin structures through weighted Poincaré-Wirtinger type inequalities. This compactness lemma, combined with the two-dimensional results of~\cite{BP1} and the cylindrical structure of the composite allows us to 
obtain an explicit formula of~${\sigma_\ast(h)}$, once again, in terms of the transversal homogenized conductivity~${\widetilde\sigma_\ast(h)}$ and of some bounded function~${\theta}$ which, in some sense, takes account of the distribution of the highly conducting phase~$\Omega_n$ in~$\Omega$~(see Theorem~\ref{thgen}).

\bigskip
The structure of the paper is the following: In Section~\ref{statement} we set up some general notations. Section~\ref{per} deals with the periodic case. In Section~\ref{gen} we extend the periodic result of Section~\ref{per} to a non periodic framework. Section~\ref{examples} is devoted to some examples illustrating both the periodic and non periodic perturbation formulas.

\medskip

Here, we give some general notations and definitions.

\subsection{General notations and definitions}   \label{statement}

\begin{itemize}
\item[$\bullet$]~${\Omega}$ is a bounded open subset of~${\R^3}$ with a Lipschitz boundary. The unit cube~${(\textstyle{-\frac{1}{2}},\frac {1}{2})^3}$ of~${\R^3}$ is denoted by~${Y}$.
\item[$\bullet$] For any subset~${\omega}$ of~${\Omega}$, we denote by~${\overline\omega}$ the closure of~${\omega}$ in~${\R^3}$.
\item[$\bullet$]~${\e_n}$ is a sequence of positive real numbers converging to zero as~${n}$ goes infinity.
\item[$\bullet$] For any matrix~${\sigma}$ in~${\R^{d \times d}}$,~${\sigma^\mathrm{T}}$ denotes the transpose of the matrix~${\sigma}$ while~${\sigma^s}$ denotes its symmetric part. For any invertible matrix~${\sigma}$ in~${\R^{d \times d}}$,~${\sigma^{-\T}:= \big( \sigma^{-1}\big)^\T = \big( \sigma^{\T}\big)^{-1}}$. 
\item[$\bullet$]~${I_d}$ denotes the unit matrix in~${\R^{d \times d}}$ and~${J:=\left(\begin{smallmatrix} 0 & -1 \\ 1 & 0 \end{smallmatrix}\right)}$.
\item[$\bullet$] For any~${h \in \R^3}$,~${\mathscr E(h)}$ denotes the~${3\times3}$ antisymmetric matrix defined by~${\mathscr E(h)\,x:=h \times x}$, for~${x\in\R^3}$.
\item[$\bullet$] For any~${\sigma,\eta\in\R^{d \times d}}$,~${\sigma\leq \eta}$ means that for any~${\xi \in \R^d}$,~${\sigma \xi \cdot \xi \leq \eta \xi \cdot \xi}$.
\item[$\bullet$] For any vector~${\xi \in \R^3}$,~${\w \xi\in\R^2}$ denotes the vector of its first two components \begin{equation*} \w \xi := (\xi_1, \xi_2)^\T.\end{equation*}
\item[$\bullet$]~${\n \cdot}$ denotes the gradient operator in~${\R^3}$ with respect to the three variables~${(x_1,x_2,x_3)}$ while~${\w \n \cdot}$ is the gradient operator in~${\R^2}$ with respect to the first two variables~${(x_1,x_2)}$: for any~${{u\in H^1(\Omega)}}$, the function~${\w \n u}$ is defined on~${\widetilde\Omega}$ by
\[
\w \n u:=\left(\frac{\partial u}{\partial x_1}, \frac{\partial u}{\partial x_2}\right)^\mathrm{T},
\]
where~${\widetilde\Omega}$ is the projection of~${\Omega}$ on the~${(x_1,x_2)}$-plane.
\item[$\bullet$] For any~${3\times3}$ matrix~${\sigma}$, we denote by~${\w \sigma}$ the~${2\times2}$ matrix defined by
\begin{equation*}
\w \sigma:=\left(\begin{smallmatrix}\sigma_{11} & \sigma_{12} \\*[.6em] \sigma_{21} & \sigma_{22} \end{smallmatrix}\right).
\end{equation*}
\item[$\bullet$] The scalar product of two vectors~${u}$ and~${v}$ of~${\R^d}$ is denoted by~${u\cdot v}$.
\item[$\bullet$]~${| \cdot |}$ denotes, the euclidean norm in~${\R^d}$, the subordinate norm in~${\R^{d \times d}}$ and the Lebesgue measure.
\item[$\bullet$] For a Borel subset~${\omega\in\R^d}$ and a function~${u\in {L}^1(\omega)}$ the average value of~${u}$ over~${\omega}$ is denoted by
\[\fint_{\omega} u \dx:=
\frac{1}{|{\omega|}}\int_{\omega}u \dx.
\] When~${\omega = Y}$, we simply denote this average value by~${\langle \cdot \rangle}$.
\item[$\bullet$] We denote by~${\mathds{1}_\omega}$ the characteristic function of the set~${\omega}$.
\item[$\bullet$] We denote by~${\mathscr C_c(\Omega)}$ the set of continuous functions with compact support in~${\Omega}$. The subspace of~${\mathscr C_c(\Omega)}$ of infinitely differentiable functions with compact support in~${\Omega}$ is denoted by~${\mathscr D(\Omega)}$.
\item[$\bullet$] We denote by~${\mathscr C_0(\Omega)}$ the space of continuous functions on~${\overline{\Omega}}$ vanishing on the boundary~${\partial\Omega}$ of~${\Omega}$ endowed with the usual norm.
\item[$\bullet$] For any locally compact subset~${X}$ of~${\R^d}$,~${\mathcal M(X)}$ denotes the set of Radon measures defined on~${X}$.
\item[$\bullet$] A sequence~${(\mu_n)}$ in~${{\cal M}(\Omega)}$ is said to weakly-${\ast}$ converge to a measure~${\mu}$ if
\begin{equation*}
\int_\Omega\varphi\mu_n(\text{d}x)\xrightarrow[n\to \infty]{}\int_\Omega\varphi\mu(\text{d}x),\quad{\rm for~any}~\varphi\in \mathscr C_0(\Omega).
\end{equation*}
\item[$\bullet$] The space of~${Y}$-periodic functions which belong to~${{L}_{\rm loc}^p(\R^d)}$~(resp.~${{H}_{\rm loc}^1(\R^d)}$) is denoted \linebreak[4] by~${{L}_{\#}^p(Y)}$~(resp.~${{H}_{\#}^1(Y)}$).
\item[$\bullet$] 
${o(\delta)}$ denotes a term of the form~${\delta\zeta(\delta)}$ where the limit of~${\mathds{\zeta}(\delta)}$ is zero, as~${\delta}$ goes to zero. For any sequences~${(a_n)_{n \in \N^*}}$ and~${(b_n)_{n \in \N^*}}$,~${a_n \underset{n \to \infty}{\sim} b_n}$ means that~${a_n = b_n + o(b_n)}$.
\item[$\bullet$] Throughout the paper, the letter~${c}$ denotes a positive constant the value of which is not given explicitly and may vary from line to line.
\end{itemize}

In the sequel, we will use the following extension of H-convergence for two-dimensional high-contrast conductivities introduced in~\cite{BrCa1} for the symmetric case and extended in~\cite{NHM} to the non symmetric case:

\begin{Def} \label{HM} 
Let~${\w \Omega}$ be a bounded domain of~${\R^2}$ and let~${\w \sigma_n \in L^\infty(\Omega)^{2 \times 2}}$ be a sequence of equi-coercive matrix-valued functions. The sequence~${\w \sigma_n}$ is said to~${H(\mathcal M(\w \Omega)^2 )}$-converge to a matrix-valued function~${\w \sigma_*}$ if for any distribution~${g}$ in~${H^{-1}(\w \Omega)}$, the solution~${u_n}$ of the problem
\begin{equation*}
\left\{\begin{array}{r l}
\divg \big(\w \sigma_n \w \nabla u_n\big)= g & \text{in } \w \Omega,
\\*[0.3em]
u_n = 0 & \text{on }\partial \w \Omega,
\end{array}\right.
\end{equation*}
satisfies the convergences
\begin{equation*}
\left\{\begin{array}{r c l l}
u_n & \harpoon & u & \text{in }H^1_0(\w \Omega),
\\*[0.3em]
\w \sigma_n \w \nabla u_n & \harpoon & \w \sigma_* \w \nabla u & \text{weakly-${*}$ in }\mathcal M(\w \Omega)^2,
\end{array}\right.
\end{equation*}
where~${u}$ is the solution of the problem
\begin{equation*}
\left\{\begin{array}{rl}
\divg \big(\w \sigma_* \w \nabla u\big)= g & \text{in }\w \Omega,
\\
u = 0 & \text{on }\partial \w \Omega.
\end{array}\right.
\end{equation*}
\end{Def}
Let~${\w \Omega}$ be a bounded open subset of~${\R^2}$ with a Lipschitz boundary and~${\w \Omega_n}$ be a sequence of open subsets of~${\w \Omega}$. Let~${\Omega}$ be the bounded open cylinder~${\Omega := \w \Omega \times (0,1)}$ and~${\Omega_n}$ the sequence of open cylinders~${\Omega_n := \w \Omega_n \times (0,1)}$.
Consider~${\alpha_1 >0}$,~${\beta_1 \in \R}$ and two sequences~${\alpha_{2,n} \geq \alpha_1}$ and~${\beta_{2,n} \in \R}$. Define, for any~${h \in \R^3}$, the two-phase isotropic conductivity \begin{equation*} 
\sigma_n(h) := \left\{ \begin{array}{l l}
\sigma_1(h):=\alpha_1 I_3 + \beta_1 \mathscr E(h)  & \text{in } \Omega \setminus \Omega_n,\\*[0.5em]
\sigma_{2,n}(h):=\alpha_{2,n} I_3 + \beta_{2,n} \mathscr E(h)  & \text{in } \Omega_n,
\end{array} \right. \quad \text{where} \quad \mathscr E(h) := \left( \begin{smallmatrix} 0 & -h_3 & h_2 \\ h_3 & 0 & -h_1 \\ -h_2 & h_1 & 0 \end{smallmatrix} \right).
\end{equation*} 
In the domain~${\Omega}$, the matrix-valued function~${\sigma_n(h)}$ does not depend on the variable~${x_3}$ and model the conductivity of a columnar heterogeneous medium. The phase~${\Omega_n}$ is the one of high conductivity:~${\alpha_{2,n}}$ and~${\beta_{2,n}}$ are unbounded. In order to ensure the~${L^1(\Omega)^{3\times3}}$-boundedness of the conductivity, we assume that the volume fraction of the highly conducting phase~${\theta_n:=|\Omega|^{-1} |\Omega_n|}$ converges to zero and that the convergences 
\begin{equation} \label{convphas20}
\left\{ \begin{array}{r c l} \theta_n \alpha_{2,n} & \xrightarrow[n \to \infty]{} & \alpha_2 >0, \\
\theta_n \beta_{2,n} & \xrightarrow[n \to \infty]{} & \beta_2 \in \R, 
\end{array}\right.
\end{equation} 
hold. Assumption~\eqref{convphas20} can be rewritten 
\begin{equation*}
\theta_n \sigma_{2,n}(h) = \theta_n \alpha_{2,n} I_3 + \theta_n \beta_{2,n} \mathscr E(h) \xrightarrow[n \to \infty]{} \sigma_2(h):= \alpha_2 I_3 + \beta_2 \mathscr E(h).
\end{equation*}
Our aim is to study the homogenization of the Dirichlet problem, for~${f \in H^{-1}(\Omega)}$,
\begin{equation} \label{Dir}
\left\{\begin{array}{r c l l}
- \divg \big( \sigma_n(h) \n u_n \big) & = & f & \text{in } \Omega, \\
u_n & = & 0 & \text{on } \partial \Omega.
\end{array}\right.
\end{equation} 

On the one hand, we consider the case of a periodic cylindrical composite without any assumption on the geometry of its cross section. This framework extends the one of the three-dimensional result of~\cite{BP1} where the highly conducting zone is a set of circular fibres. On the other hand, by the means of a compactness result~(see Lemma~\ref{divcurl}), we analyse the case of cylindrical but non periodic composites. In both cases, we impose conditions, adapting~\cite{BrARMA}, that prevent from the appearance of non local terms so that the limit equation of~\eqref{Dir} is a conductivity one. 

\medskip

In the sequel, we will omit the dependence on~$h$ of~${\sigma_1(h)}$,~${\sigma_{2,n}(h)}$ and~${\sigma_2(h)}$ denoting simply~${\sigma_1}$,~${\sigma_{2,n}}$ and~${\sigma_2}$. 


 
\markboth{Chapitre~\ref*{chap2} - Homogenization of high-contrast and non symmetric conductivities for non periodic columnar structures}{\textit{Section~\ref*{per} - The periodic case}}

\section{The periodic case} \label{per}

\markboth{Chapitre~\ref*{chap2} - Homogenization of high-contrast and non symmetric conductivities for non periodic columnar structures}{\textit{Section~\ref*{per} - The periodic case}}

In this section, we study the influence of a constant magnetic field~${h \in \R^3}$ on the effective conductivity of a composite material where the highly conducting phase is periodically distributed 
 but, contrary to~\cite{BP1}, the cross section of which has a general geometry. Consider a sequence~${\omega_n = \w \omega_n \times (0,1)}$ where~${\w \omega_n}$ is a sequence of subsets of~${(0,1)^2}$ with~${|\omega_n|}$ converging to 0, as~${n}$ tends to infinity. Let~${\Omega_n}$ be the sequence of open subsets of~${\Omega}$ defined by 
\begin{equation*} 
\Omega_n = \Omega \cap \bigcup_{k \in \mathbb Z^3} \e_n \big( \omega_n + k\big).
\end{equation*} 
The conductivity of the heterogeneous medium occupying~${\Omega}$ is given by 
\begin{equation} \label{defper}
\sigma_n(h)(x) = \Sigma_n(h)\left( \frac{x}{\e_n}\right), \quad \forall x \in \Omega,
\end{equation}
where~${\Sigma_n(h)(\cdot)}$ is a~${Y}$-periodic function defined by 
\begin{equation} \label{defperan}
\Sigma_n(h)= a_n I_3 + b_n \mathscr E(h) \quad \text{with} \quad \left\{ \begin{array}{l l} a_n:= \alpha_1 \mathds{1}_{Y \setminus \omega_n} + \alpha_{2,n}\mathds{1}_{\omega_n}, \\ 
b_n:= \beta_1\mathds{1}_{Y \setminus \omega_n} + \beta_{2,n}\mathds{1}_{\omega_n}.
\end{array} \right.
\end{equation} 

For a fixed~${n \in \N^*}$, let~${(\sigma_n)_*(h)}$ be the constant matrix defined by
\begin{equation} \label{defsigman*}
\forall\,\lambda \in \R^3, \quad (\sigma_n)_*(h) \lambda =\big\langle \Sigma_n(h) \nabla W_n^\lambda \big\rangle, 
\end{equation}
where, for any~${\lambda \in \R^3}$,~${W_n^\lambda}$ is the unique solution in~${ H^1_\sharp(Y)}$ of the auxiliary problem
\begin{equation} \label{pb3}
\displaystyle \divg \big( \Sigma_n(h) \nabla W_n^\lambda\big) = 0 \quad \text{in } \mathscr{D}'(\R^3) \quad \text{and} \quad \big\langle W_n^\lambda - \lambda \cdot y \big\rangle =0,
\end{equation}
which is equivalent to the variational cell problem
\begin{equation} \label{pb4}
\left\{\begin{array}{r c l l}
\displaystyle &\displaystyle\big\langle \Sigma_n(h) \nabla W_n^\lambda \cdot \nabla \Phi \big\rangle = 0, \quad \forall\,\Phi \in H^1_{\sharp}(Y), \vspace{0.2cm}\\
&\displaystyle\big\langle W_n^\lambda(y) - \lambda \cdot y \big\rangle = 0.
\end{array} \right.
\end{equation} 
The matrix~${(\sigma_n)_*(h)}$ is the homogenized conductivity of the oscillating sequence~${\Sigma_n(\cdot / \e)}$ as~${\e \to 0}$~(see, for instance,~\cite{BLP} for more details).


The limit problem of the high-contrast three-dimensional equation~\eqref{Dir} where~${\sigma_n(h)}$ is given by~\eqref{defper} may include non local effects. In order to avoid such effects, we assume, following~\cite{BrARMA}, that the weighted Poincaré-Wirtinger inequality
\begin{equation} \label{WPWper}
\forall V \in H^1(Y), \quad \int_{Y} a_n \left| V - \int_Y V\right|^2 \leq C_n \int_{Y} a_n |\nabla V|^2, \end{equation} holds true with 
\begin{equation}  \label{condWPWper}
\e_n^2 C_n \xrightarrow[n \to \infty]{} 0. 
\end{equation} 
Under the assumptions~\eqref{WPWper} and~\eqref{condWPWper}, it was shown in~\cite{BrARMA} that the sequence of problems~\eqref{Dir} converges to a conduction one with a homogenized conductivity~${\sigma_*(h)}$. 

\medskip

The main contribution of Proposition~\ref{propper} below is to provide a formula for the effective conductivity~${\sigma_*(h)}$ 
of a cylindrical periodic composite the cross section of which has a general geometry.



\begin{Prop} \label{propper}
Consider the sequence of problems~\eqref{Dir} where~${\sigma_n(h)}$ is the conductivity defined by~\eqref{defper}-\eqref{defperan}. Assume that~\eqref{convphas20},~\eqref{WPWper} and~\eqref{condWPWper} are satisfied. Then, there exists a constant matrix~${\sigma_*(h)}$ such that, up to a subsequence, the solution~${u_n}$ of~\eqref{Dir} weakly converges in~${H_0^1(\Omega)}$ to the solution~${u}$ of
\begin{equation} \label{Dir*per}
\left\{\begin{array}{r c l l}
- \divg \big( \sigma_*(h) \n u \big) & = & f & \text{in } \Omega, \\
u & = & 0 & \text{on } \partial \Omega.
\end{array}\right.
\end{equation} 
Moreover, the homogenized matrix~${\sigma_*(h)}$ is the limit of~${\big(\sigma_n\big)_*(h)}$~(see~\eqref{defsigman*}) and is given by 
\begin{equation} \label{propper'}
\sigma_*(h) := \begin{pmatrix}
\w \sigma_* & p_* \\*[0.2em]
q_*^\T & \alpha_*
\end{pmatrix},
\end{equation} 
where
\begin{equation} \label{eqper}
\left\{ \begin{array}{l}
\vspace{0.2cm} \displaystyle p_* = - \left[\beta_1 I_2 + \beta_2 \big(\w \sigma_* - \w \sigma_1 \big) \w \sigma_2^{-1}\right] J \w h, \\
\vspace{0.2cm} \displaystyle q_* = \left[\beta_1 I_2 + \beta_2 \, \w \sigma_2^{-1} \big(\w \sigma_* - \w \sigma_1 \big)\right]^\T J \w h,\\
\alpha_* = \alpha_1 + \alpha_2 +  \beta_2^2 \, \w \sigma_2^{-1} \big(\w \sigma_1 + \w \sigma_2  - \w \sigma_*\big) \w \sigma_2^{-1} J \w h \cdot J \w h,
\end{array}\right.
\end{equation}
and, for any~${i=1,2}$,
\begin{equation*}
\w \sigma_i := \begin{pmatrix} \alpha_i & - \beta_i \, h_3 \\*[0.3em]
\beta_i \, h_3 & \alpha_i \end{pmatrix}.
\end{equation*}
\end{Prop}

\begin{Rem} \upshape
For the sake of simplicity, throughout the paper, the symmetric part of~${\sigma_n(h)}$ is supposed to be isotropic. However, the results we obtain can be extended to composites the components of which have anisotropic conductivities.
\end{Rem}

\begin{Rem} \upshape \label{remper}
It was shown in~\cite{BrARMA} that, due to the~${L^1(Y)^{3 \times 3}}$-boundedness of~${\Sigma_n(h)(\cdot)}$, the sequence~${\big(\sigma_n\big)_*(h)}$ is bounded. Thanks to~\eqref{WPWper} and~\eqref{condWPWper}, Theorem 2.1 of~\cite{BrARMA} ensures that the limit~${\sigma_*(h)}$ obtained in the following way
\begin{equation*}
\Sigma_n(h)\left(\frac{x}{\e}\right) \underset{\e \to 0}{\overset{H}{\harpoon}} \big(\sigma_n\big)_*(h) \xrightarrow[n \to \infty]{} \sigma_*(h),
\end{equation*}
satisfies the convergence
\begin{equation*} 
\sigma_n(h) \n u_n \harpoon \sigma_*(h) \n u \quad \text{in } \mathcal M(\Omega)^3,
\end{equation*}
and then, coincides with the homogenized conductivity matrix in the problem~\eqref{Dir*per}.
\end{Rem}

\begin{Rem} \upshape \label{remexpl}
Since~${\Omega_n}$ has a columnar structure, the sequence~${\w \sigma_n(h)}$ given by \begin{equation*}
\w \sigma_n(h) := \w \sigma_n(h_3) = \left\{ \begin{array}{l l}
\w \sigma_1(h_3)=\alpha_1 I_2 + \beta_1 h_3 J   & \text{in } \w \Omega \setminus \w \Omega_n,\\*[.3em]
\w \sigma_{2,n}(h_3)=\alpha_{2,n} I_2 + \beta_{2,n} h_3 J  & \text{in } \w \Omega_n.
\end{array} \right.
\end{equation*}
depends only on the transversal variable~${(x_1,x_2)}$ and is then associated with the two-dimensional problems, for any~${g \in H^{-1}(\w \Omega)}$,
\begin{equation*}
\left\{\begin{array}{r c l l}
- \divg \big( \w \sigma_n(h_3) \w \n v_n \big) & = & g & \text{in } \w \Omega, \\*[.3em]
v_n & = & 0 & \text{on } \partial \w \Omega.
\end{array}\right.
\end{equation*} 
Similarly to~\eqref{defsigman*}, we define the constant matrix~${\big(\w \sigma_n\big)_*(h_3)}$. For any~${\lambda \perp e_3}$, the solution~${W_n^\lambda}$ of~\eqref{pb3} does not depend on the variable~${y_3}$ and then 
\begin{equation*}
\displaystyle \big\langle W_n^\lambda - \w \lambda \cdot \w y \big\rangle =0 \quad \text{and} \quad \divg \big( \w \Sigma_n(h_3) \w \nabla W_n^\lambda\big) = 0 \quad \text{in} \ \mathscr{D}'(\R^2).
\end{equation*} 
This equation and~\eqref{pb3} imply that, for any~${\lambda, \mu \perp e_3}$, 
\begin{equation*} 
\big(\sigma_n\big)_*(h)\lambda \cdot \mu = \big\langle \Sigma_n(h) \n W_n^\lambda \big\rangle \cdot \mu  = \big\langle \w \Sigma_n(h_3) \w \n W_n^\lambda \big\rangle \cdot \w \mu = \big( \w \sigma_n\big)_*(h_3)\w \lambda \cdot \w \mu. 
\end{equation*} 
Hence, by Remark~\ref{remper},~${\big( \w \sigma_n\big)_*(h_3)}$ converges to the~${2 \times 2}$ matrix~${\w \sigma_*}$ involved in~\eqref{propper'}. 
A two-dimensional perturbation formula in~\cite{BP1} gives the influence of the magnetic field~${h_3}$ on~${\w \sigma_*}$: 
\begin{equation} \label{wsigma}
 \displaystyle \w \sigma_*:= \w \sigma_*(h_3) =\sigma^0_*\big(\alpha_1, \alpha_2 + \alpha_2^{-1} \beta_2^2 h_3^2 \big) + h_3 \beta_1 J,
\end{equation} 
where~${\sigma_*^0}$ is a locally Lipschitz function defined on~${(0, \infty)^2}$, and for any~${\alpha_1, \alpha_2 >0}$,~${\sigma_*^0(\alpha_1, \alpha_2)}$ is the transversal homogenized conductivity in the absence of a magnetic field. The independence of the microstructure of the variable~${x_3}$ allows us to obtain an explicit expression of~${\sigma_*(h)}$ in terms of the transversal homogenized conductivity~${\sigma_*^0}$ 
in the absence of a magnetic field.
\end{Rem}

\begin{Rem} \upshape
In the case where the high conducting phase is a set of circular fibres, it was proved in~\cite{BP1} that~${\sigma_*^0(\alpha_1,\alpha_2) = \alpha_1 I_2}$ and the limit~${\sigma_*(h)}$ in~\eqref{propper'} reduces to 
\begin{equation*} 
\sigma_*(h) = \alpha_1 I_3 + \left[ \alpha_2 +  \beta_2^2 \, \w \sigma_2^{-1} J \w h \cdot J \w h \right] e_3 \otimes e_3 + \beta_1 \mathscr E(h).
\end{equation*}
\end{Rem} 

Now, let us proceed with the proof of Proposition~\ref{propper}.

\medskip

\noindent \textbf{Proof of Proposition~\ref{propper}.} Thanks to Remarks~\ref{remper} and~\ref{remexpl}, there exists a~${{3 \times 3}}$ matrix~${\sigma_*(h)}$ such that, up to a subsequence, we have the convergence of constant matrices 
\begin{equation} \label{conv0}
\big( \sigma_n\big)_*(h) \xrightarrow[n \to \infty]{}  \sigma_*(h) := \begin{pmatrix}
\w \sigma_* & p_* \\*[.2em]
q_*^\T & \alpha_*
\end{pmatrix},
\end{equation} 
where~${\w \sigma_*}$ is given by~\eqref{wsigma} and where the constants~${q_*, p_*} \in \R^2$,~${\alpha_*} \in \R$ have to be determined. To this end, we divide the proof into two steps. We first apply Grabovsky and Milton's method~\mbox{\cite{Grab1,Grab2,MilGrab}} to link~${\big( \sigma_n\big)_*}$ to a more simple problem. Then, we study the asymptotic behavior of the different coefficients of this new problem.

\bigskip

\noindent \textit{First step:} A stable transformation under homogenization. For a fixed~${n \in \N^*}$, following Grabovsky and Milton~\mbox{\cite{MilGrab,Grab1}}, we consider two vectors~${p_{0,n}, q_{0,n} \in \R^2}$ and the transformation 
\begin{equation} \label{k1}
\sigma'_n := {\Pi}_n \ \sigma_n(h) \ \widehat \Pi_n = \begin{pmatrix}
\w \sigma_n & p'_n \\*[.4em]
{q'}_n^\T & \alpha'_n
\end{pmatrix},
\end{equation} 
where 
\begin{equation}  \label{defPi} 
{\Pi}_n:=\begin{pmatrix}
I_2 & 0 \\*[.4em]
q_{0,n}^\T & 1 
\end{pmatrix}, \quad \widehat \Pi_n:=\begin{pmatrix}
I_2 & p_{0,n} \\*[.4em]
0  & 1 
\end{pmatrix},
\end{equation} 
and 
\begin{equation}  \label{defpn'}
p'_n = \left\{ \begin{array}{l l}
\w \sigma_1 p_{0,n} - \beta_1 J \w h & \text{in } \Omega \setminus \Omega_n, \\
\w \sigma_{2,n} p_{0,n} - \beta_{2,n} J \w h & \text{in } \Omega_n,
\end{array}\right. \quad q'_n = \left\{ \begin{array}{l l}
\w \sigma_1^\T q_{0,n} + \beta_1 J \w h & \text{in } \Omega \setminus \Omega_n, \\
\w \sigma_{2,n}^\T q_{0,n} + \beta_{2,n} J \w h & \text{in } \Omega_n.
\end{array}\right. 
\end{equation} 
Let us choose the parameters~${p_{0,n}}$ and~${q_{0,n}}$ in such a way that~${p'_n}$ and~${q'_n}$ are constant. To that aim,~${p_{0,n}}$ and~$q_{0,n}$ have to satisfy the identities 
\begin{equation*} 
\left\{ \begin{array}{r c l} \w \sigma_1 p_{0,n} - \beta_1 J \w h & = & \w \sigma_{2,n} p_{0,n} - \beta_{2,n} J \w h, \\*[0.5em] 
\w \sigma_1^\T q_{0,n} + \beta_1 J \w h & = & \w \sigma_{2,n}^\T q_{0,n} + \beta_{2,n} J \w h,
\end{array} \right.\end{equation*} 
which implies that 
\begin{equation*}
p_{0,n} = (\beta_{2,n} - \beta_1)\big( \w \sigma_{2,n} - \w \sigma_1 \big)^{-1} J \w h \quad \text{and} \quad q_{0,n} = (\beta_{2,n} - \beta_1)\big( \w \sigma_1  - \w \sigma_{2,n}\big)^{-\T} J \w h.
\end{equation*} 
The new matrix-valued function~${\sigma'_n}$ defined by~\eqref{k1} is periodic and can be rewritten
\begin{equation} \label{alpha''''}
\forall x \in \Omega, \quad \sigma'_n(x) = \Sigma'_n\left( \frac{x}{\e_n}\right) \quad \text{where} \quad \Sigma'_n:= \begin{pmatrix}
\w \Sigma_n & p'_n \\*[.4em]
{q'}_n ^\T & a'_{n}
\end{pmatrix}\! . 
\end{equation} 
Moreover, by~\eqref{k1}, the coefficient~${a'_n}$ in~\eqref{alpha''''} has the following explicit expression: 
\begin{equation} \label{alpha'} 
a'_n =  \alpha'_{1,n} \mathds{1}_{Y \setminus \omega_n} + \alpha'_{2,n} \mathds{1}_{\omega_n} \quad \text{where} \quad \left\{\begin{array}{l} \vspace{0.2cm}\alpha'_{1,n} = \alpha_1 + \w \sigma_1 p_{0,n} \cdot q_{0,n} + \beta_1 \, (p_{0,n} - q_{0,n}) \cdot J \w h, \\*[.4em]
\alpha'_{2,n} = \alpha_{2,n} + \w \sigma_{2,n} p_{0,n} \cdot q_{0,n} + \beta_{2,n} \, (p_{0,n} - q_{0,n}) \cdot J \w h.\end{array}\right.
\end{equation} 

Let us now study the homogenization of~${\sigma'_n}$. Define~${\big(\sigma'_n\big)_*}$ as in the formula~\eqref{defsigman*}. The conductivity~${\Sigma'_n}$ does not depend on the variable~${y_3}$. On the one hand, as in Remark~\ref{remexpl}, if~${\lambda \perp e_3}$, the solution~${W_n^\lambda}$ of the problem~\eqref{pb4}, with the conductivity~${\Sigma'_n}$, does not depend on the variable~${y_3}$ and~${\nabla W_n^\lambda = ( \w \n W_n^\lambda,0)^\T}$. Hence, since~${q'_n}$ is a constant, and by Remark~\ref{remexpl}, 
\begin{equation} \label{sigm'2}
\big(\sigma'_n \big)_* \lambda = \big( \langle \w \Sigma_n \w \n W_n^\lambda \rangle, \langle q'_n \cdot \w \n W_n^\lambda \rangle \big)^\T = \Big(\big(\w \sigma_n\big)_* \w \lambda, q'_n \cdot \w \lambda\Big)^\T.
\end{equation}  
On the other hand, it is clear that, for~${\lambda = e_3}$,~${W_n^{e_3}(y) = y_3}$ satisfies~\eqref{pb4} with the conductivity~${\Sigma'_n}$. Hence, since~${p'_n}$ is a constant, we have 
\begin{equation} \label{sigm'1}
\big(\sigma'_n \big)_* e_3 = \big\langle (p'_n, a'_n)^\T \big\rangle = \big( p'_n, \langle a'_n\rangle\big)^\T.
\end{equation} 
Then, by~\eqref{alpha'},~\eqref{sigm'2},~\eqref{sigm'1} and since~${|\omega_n| \underset{n \to \infty}{\sim} \theta_n}$, the matrix~${\big(\sigma'_n \big)_*}$ has the form 
\begin{equation} \label{sigman'}
\big(\sigma'_n \big)_* = \begin{pmatrix}
\vspace{0.1cm}\big(\w \sigma_n \big)_* & p'_n \\
{q'}_n ^\T & \langle a'_{n} \rangle
\end{pmatrix}, 
\end{equation} 
where 
\begin{equation} \label{alpha''}
\begin{array}{r c l} \langle a'_{n} \rangle & = & \left[ \alpha_1 + \w \sigma_1 p_{0,n} \cdot q_{0,n} + \beta_1 \, (p_{0,n} - q_{0,n}) \cdot J \w h \right] \\*[0.7em] 
& + & \theta_n \left[ \alpha_{2,n} + \w \sigma_{2,n} p_{0,n} \cdot q_{0,n} + \beta_{2,n} \, ( p_{0,n} - q_{0,n}) \cdot J \w h  \right] + o(1). \end{array}
\end{equation}

\bigskip

\noindent \textit{Second step:} Application of the theory of exact relations and asymptotic behavior of~${\big(\sigma'_n \big)_*}$. By~\eqref{convphas20} and since the volume fraction~${\theta_n}$ converges to 0, we have 
\begin{equation} \label{conv1'}
\left\{ \begin{array}{r c l} p_{0,n} = \theta_n (\beta_{2,n} - \beta_1)\big(\theta_n (\w \sigma_{2,n} - \w \sigma_1) \big)^{-1} J \w h & \xrightarrow[n \to \infty]{} & \beta_2 \w \sigma_2^{-1} J \w h \\*[0.7em]
q_{0,n}  = \theta_n (\beta_{2,n} - \beta_1)\big( \theta_n (\w \sigma_1  - \w \sigma_{2,n})\big)^{-\T} J \w h & \xrightarrow[n \to \infty]{} & - \beta_2 \w \sigma_2^{-\T} J \w h.
\end{array} \right. 
\end{equation} 
Then, by~\eqref{convphas20},~\eqref{defpn'},~\eqref{defPi},~\eqref{alpha''} and~\eqref{conv1'} we obtain the convergences 
\begin{equation} \label{conv1}
\left \{\begin{array}{r c l} p'_n & \underset{n \to \infty}{\longrightarrow}  & p'_* := \left[- \beta_1 I_2 + \beta_2 \w \sigma_1 \w \sigma_2^{-1}\right] J \w h, \\*[0.7em]
q'_n & \underset{n \to \infty}{\longrightarrow}  & q'_* := \left[\beta_1 I_2 - \beta_2 \w \sigma_1^\T \w \sigma_2^{-\T}\right] J \w h, \\*[0.9em]
 \langle a'_{n} \rangle & \underset{n \to \infty}{\longrightarrow} & \alpha'_* := \displaystyle \sum_{i=1}^2 \left[\alpha_i - \beta_2^2 \, \w \sigma_2^{-1} \w \sigma_i \w \sigma_2^{-1} J \w h \cdot J \w h + 2 \, \beta_2 \, \beta_i \, \w \sigma_2^{-1} J \w h \cdot J \w h \right], \end{array} \right. 
\end{equation} 
\begin{equation} \label{zz}
{\Pi}_n\xrightarrow[n \to \infty]{} \Pi:=\begin{pmatrix}
I_2 & 0\\
\beta_2 \w h^\T J \w \sigma_2^{-1}  & 1 
\end{pmatrix} \qquad \text{and} \qquad \widehat \Pi_n\xrightarrow[n \to \infty]{} \widehat \Pi:=\begin{pmatrix}
I_2 & \beta_2 \w \sigma_2^{-1} J \w h  \\
0  & 1 
\end{pmatrix}. 
\end{equation}

Since the matrix transformation~\eqref{k1} preserves the H-limit in the periodic case~(see, for instance,~\mbox{\cite{Grab1,Grab2,MilBook}}), we have 
\begin{equation} \label{ident'}
\big(\sigma'_n\big)_* = {\Pi}_n \ \! \big(\sigma_n\big)_*(h) \ \widehat \Pi_n. 
\end{equation} 
Passing to the limit, as~${n}$ goes to infinity, in relation~\eqref{ident'}, using~\eqref{conv0},~\eqref{sigman'},~\eqref{conv1}-\eqref{zz}, we obtain 
\begin{equation} \label{ident}
\begin{pmatrix}
\w \sigma_* & p'_* \\
{q'}_*^\T & \alpha'_* 
\end{pmatrix} =   \Pi \ \sigma_* \ \widehat \Pi.
\end{equation} 
Inverting the identity~\eqref{ident} and taking into account~\eqref{conv1} and~\eqref{zz} leads to~\eqref{eqper}. The proof of Proposition~\ref{propper} is completed. 
\qed

\medskip

Now let us turn to the non periodic case.

\markboth{Chapitre~\ref*{chap2} - Homogenization of high-contrast and non symmetric conductivities for non periodic columnar structures}{\textit{Section~\ref*{gen} - The non periodic case}}


\section{The non periodic case}\label{gen}

\markboth{Chapitre~\ref*{chap2} - Homogenization of high-contrast and non symmetric conductivities for non periodic columnar structures}{\textit{Section~\ref*{gen} - The non periodic case}}

In this section, we study the homogenization of the problem~\eqref{Dir} without any periodicity assumption. The conductivity~${\sigma_n(h)}$ is defined by 
\begin{equation} \label{defgeneq}
\sigma_n(h) := \alpha_n I_3 + \beta_n \mathscr E(h) \quad \text{where} \quad \left\{ \begin{array}{l l} \alpha_n:= \1 \alpha_1 + \2 \alpha_{2,n}, \\ 
\beta_n:= \1 \beta_1 + \2 \beta_{2,n}.
\end{array} \right.
\end{equation}  
Consider the covering of~${\R^3}$ by the squares~${Q_n^k}$ defined by
\begin{equation} \label{decoup}
\forall k \in \mathbb Z^3, \quad Q_n^k = \e_n (Y + k).
\end{equation}
We assume that the conductivity coefficient~${\alpha_n}$ defined by~\eqref{defgeneq} satisfies, for any~${k \in \Z^3}$,~${n\in\N^*}$, the following conditions: 
\begin{enumerate}
\item[$(i)$] the weighted Poincaré-Wirtinger inequality 
\begin{equation} \label{WPW}
 \forall \, v \in H^1(Q_n^k), \quad \int_{Q_n^k} \alpha_n \left| v - \fint_{Q_n^k} v\right|^2 \dx \leq c_n \int_{Q_n^k} \alpha_n |\nabla v|^2 \dx,
\end{equation} where~${c_n}$ is a sequence of positive constants satisfying 
\begin{equation} \label{condWPW} 
c_n \xrightarrow[n \to \infty]{} 0 ; 
\end{equation}
\item[$(ii)$] there exists a positive constant~${c}$ such that, for any~${k \in \Z^3}$ and~${n \in \N^*}$, 
\begin{equation} \label{ii}
\fint_{Q_n^k} \alpha_n \leq c.
\end{equation}
\end{enumerate} 

\begin{Rem} \upshape \label{noncon}
Note that, in the periodic case, the hypothesis~\eqref{WPWper}-\eqref{condWPWper} is a rescaling of \mbox{\eqref{WPW}-\eqref{condWPW}} which, similarly to the periodic case, 
prevents from the appearance of non local effects in the limit problem.
Assumption~\eqref{ii} ensures that the microstructure does not concentrate on a lower dimension subset through the homogenization process since it implies that~(see in the proof of Lemma~\ref{divcurl})
\begin{equation} \label{tt}
\theta_n^{-1} \2 \harpoon \theta \in L^\infty(\Omega) \quad \text{weakly-${*}$ in } \mathcal M(\Omega).
\end{equation} 
In the periodic case,~\eqref{ii} is clearly satisfied since 
\begin{equation*} 
\fint_{Q_n^k} \alpha_n \dx = \Vert a_n\Vert_{L^1(Y)} \leq c,
\end{equation*} 
where~${a_n}$ is defined by~\eqref{defperan} and~${\theta \equiv 1}$. 
\end{Rem}

\noindent We have the following result:

\begin{Theo} 
Assume that~\eqref{convphas20},~\mbox{\eqref{decoup}-\eqref{ii}} are satisfied. Then, there exist a matrix-valued function~${\sigma_*(h)}$ and a subsequence of~${n}$, still denoted by~${n}$, such that the solution~${u_n}$ of the problem~\eqref{Dir} converges weakly in~${H_0^1(\Omega)}$ to the solution~${u}$ of the conductivity problem
\begin{equation} \label{eq*gen}
\left\{\!\! \begin{array}{r c l l}
- \divg \big(\sigma_*(h) \nabla u \big) & = & f &\text{in } \Omega,\\
u & = & 0 & \text{on } \partial \Omega.
\end{array} \right. \;\;\; 
\end{equation} 
Moreover, 
the effective conductivity~${\sigma_*(h)}$ in~\eqref{eq*gen} is given by 
\begin{equation} \label{propgen'}
\sigma_*(h) := \begin{pmatrix}
\w \sigma_* & p_* \\*[0.2em]
q_*^\T & \alpha_*
\end{pmatrix}, 
\end{equation} 
where~${\w \sigma_*}$ is the~${H(\mathcal M(\w \Omega)^2)}$-limit of~${\w \sigma_n(h)}$ in the sense of Definition~\ref{HM},~${\theta \in L^\infty(\Omega)}$ is the weak-${*}$ limit of~${\theta_n^{-1} \2}$ and 
\begin{equation} \label{eqgen}
\left\{ \begin{array}{l}
\vspace{0.2cm} \displaystyle p_* = - \left[\beta_1 I_2 + \beta_2 \big(\w \sigma_* - \w \sigma_1 \big) \w \sigma_2^{-1}\right] J \w h, \\
\vspace{0.2cm} \displaystyle q_* = \left[\beta_1 I_2 + \beta_2 \, \w \sigma_2^{-1} \big(\w \sigma_* - \w \sigma_1 \big) \right]^\T J \w h, \\
\alpha_* = \alpha_1 + \theta \alpha_2 +  \beta_2^2  \ \! \w \sigma_2^{-1} \big(\w \sigma_1 + \theta \w \sigma_2  - \w \sigma_*\big) \w \sigma_2^{-1} J \w h \cdot J \w h,
\end{array}\right.
\end{equation} 
and, for any~${i=1,2}$,  
\begin{equation*}
\w \sigma_i := \begin{pmatrix} \alpha_i & - \beta_i \, h_3 \\*[0.3em]
\beta_i \, h_3 & \alpha_i \end{pmatrix}.
\end{equation*}
\label{thgen}
\end{Theo}

\begin{Rem} \upshape The shape~\eqref{decoup} of~${Q_n^k}$ is purely technical and can be generalized into any subset the diameter of which is of order~${\e_n}$. 
\end{Rem}

\begin{Rem} \upshape \label{remexpl'} 
Since~${\Omega_n}$ has a columnar structure,~${\mathds{1}_{\Omega_n}}$ does not depend on the variable~${x_3}$. Therefore, 
\begin{equation} \label{tt'}
\theta_n^{-1} \mathds{1}_{\w \Omega_n} \harpoon \theta \in L^\infty(\w \Omega) \quad \text{weakly-${*}$ in } \mathcal M(\w \Omega).
\end{equation} 
Hence, as in Remark~\ref{remexpl}, it was proved in~\cite{BP1} that there exists a function~${\sigma_*^0}$ defined on~${(0,\infty)^2}$ and a subsequence of~${n}$, such that, for any~${\alpha_1, \alpha_2 >0}$,~${\beta_1, \beta_2 \in \R}$, 
\begin{equation*} 
\w \sigma_n(h) = \w \sigma_n(h_3) \overset{H(\mathcal M (\w \Omega)^2)}{\harpoon} \w \sigma_*(h)=\sigma^0_*\big(\alpha_1, \alpha_2 + \alpha_2^{-1} \beta_2^2 h_3^2 \big) + h_3 \beta_1 J.
\end{equation*}
We obtain, once again, an explicit expression of~${\sigma_*(h)}$ in terms of the homogenized perturbed conductivity in the~${(x_1,x_2)}$-plane, in the absence of a magnetic field.
\end{Rem}

\noindent A crucial ingredient of the proof of Theorem~\ref{thgen} is the following compactness result:

\begin{Lem}  \label{divcurl}
Let~${\alpha_n}$ be the sequence defined by~\eqref{defgeneq} such that~\eqref{convphas20} and~\eqref{WPW}-\eqref{ii} hold true. Consider two sequences~${\xi_n \in L^1(\Omega)}$ and~${v_n \in H^1(\Omega)}$ satisfying 
\begin{equation} \label{convlem}
\xi_n \harpoon \xi \text{ weakly-${*}$ in } \mathcal{M}(\Omega) \quad \textrm{and} \quad v_n \harpoon v \text{ weakly in } H^1(\Omega).
\end{equation}
We assume that
\begin{equation} \label{divrot1}
\int_\Omega \alpha_n^{-1} |\xi_n|^2 \dx + \int_\Omega \alpha_n |\nabla v_n|^2 \dx \leq c. 
\end{equation} 
Then,~${\xi \in L^2(\Omega)}$ and we have the convergence, in the sense of distributions 
\begin{equation} \label{convdivcurl}
\xi_n  v_n \harpoon \xi v \quad \text{in } \mathscr{D}'(\Omega). 
\end{equation}
\end{Lem}

\begin{Rem} \upshape
Note that Lemma~\ref{divcurl} is false when the conditions~\eqref{WPW} and~\eqref{condWPW} do not hold. 
This can be seen by considering the classical model example of non local effects in conduction due to Fenchenko-Khruslov~\cite{FeKh} and presented, for instance, in~\mbox{\cite{BrCa1,Brinonloc}}. For the reader's convenience, we give the main steps of the counterexample. Let~${\Omega:=(\textstyle{-\frac 1 2},\textstyle{\frac 1 2})^3}$ and~$\Omega_n$ be the~\mbox{${\scriptscriptstyle{ \frac{1}{n}}}$-periodic} lattice of thin vertical cylinders of radius~${n^{-1} e^{-n^2}}$. Let~$\alpha_n$ be the conductivity defined by~\eqref{defgeneq} with~${\alpha_1 := 1}$ and~${\alpha_{2,n} := \pi^{-1} \, e^{2 n^2}}$ which satisfies~\eqref{convphas20} and~\eqref{ii}. For a fixed~$f$ in~$L^2(\Omega)$, let~$u_n$ be the solution, in~${H_0^1(\Omega)}$, of the equation 
\begin{equation*}
- \divg ( \alpha_n \n u_n ) = f \quad \text{in } \mathscr D'(\Omega).
\end{equation*}
For~${R \in (0,\textstyle{ \frac 1 2 })}$, let~${ V_n}$ be the~\mbox{$Y$-periodic} function defined on~${\R^3}$ by
\begin{equation*}
 V_n(y):= \left\{ \begin{array}{c l}
\displaystyle \frac{\ln r + n^2}{\ln R + n^2} & \text{if} \quad r:= \sqrt{y_1^2 +y_2^2} \in (e^{-n^2},R), \\*[0.7em]
0 & \text{if} \quad r \leq e^{-n^2} \text{ (region of high conductivity)}, \\*[0.7em]
1 & \text{if} \quad r \geq R.
\end{array} \right.
\end{equation*} 
An easy computation shows that the sequences~${\xi_n:=\alpha_n \n u_n \cdot e_3}$ and~${v_n(x):= V_n(nx)}$ satisfy the assumption~\mbox{\eqref{divrot1}} and that~${v_n}$ weakly converges to the constant function~${1}$ in~${H^1(\Omega)}$. Moreover, Briane and Tchou~\cite{BrTc} proved that
\begin{equation*} 
\xi_n = \alpha_n \displaystyle \frac{\partial u_n}{\partial x_3} \harpoon \displaystyle  \xi := \frac{\partial u}{\partial x_3} + \frac{\partial v}{\partial x_3} \quad \text{weakly-$*$ in } \mathcal M(\Omega), 
\end{equation*} 
where the weak limit~$u$ of~$u_n$ in~${H_0^1(\Omega)}$ and the weak-$*$ limit~${v}$ of~${{\frac{\2}{\pi e^{-2 n^2}}} \, u_n}$ in the sense of Radon measures satisfy the coupled system 
\begin{equation} \label{coupled}
\left\{ \begin{array}{r c l l}
-  \Delta u + 2 \pi \, (u-v) & = & f & \text{in} \quad \Omega, \\*[0.3em]
- \displaystyle \frac{\partial^2 v}{\partial x_3^2} + v-u & = & 0 & \text{in} \quad \Omega, \\*[0.7em]
 u & = & 0 & \text{on} \quad \partial \Omega, \\*[0.7em]
v(x',0) = v(x',1) & = & 0 & \text{if} \quad x' \in (\textstyle{-\frac 1 2},\textstyle{\frac 1 2})^2. 
\end{array} \right. 
\end{equation}
Then, if~${f}$ is non zero,~${u}$ and~${v}$ are two different functions. Therefore, the convergence~\eqref{convdivcurl} does not hold true since, by the strong convergence, up to a subsequence, of~${v_n}$ to~$1$ in~${L^2(\Omega)}$ and the weak convergence of~${\1 \n u_n}$ to~$\n u$ in~${L^2(\Omega)}$, we have  
\begin{equation*}
\xi_n v_n = \1 \displaystyle \frac{\partial u_n}{\partial x_3} \, v_n \harpoon \displaystyle \frac{\partial u}{\partial x_3} \neq \xi \times 1 = \frac{\partial u}{\partial x_3} + \frac{\partial v}{\partial x_3} \quad \text{in } \mathscr D'(\Omega).
\end{equation*} 
Substituting the expression of~${v}$, in terms of~${u}$, in the first equation of~\eqref{coupled} leads to a non local term in the equation satisfied by~${u}$. The Poincaré-Wirtinger control~\mbox{\eqref{WPW}-\eqref{condWPW}} is fundamental to avoid such effects. In this example,~\eqref{condWPW} is false since~(see~\cite{BrARMA} for more details) the optimal constant~${c_n}$ in~\eqref{WPW} satisfies~${c_n \geq c >0}$.   
\end{Rem} 

\noindent \textbf{Proof of Lemma~\ref{divcurl}.} On the one hand, by~\eqref{convphas20}, the sequence~${\alpha_n}$ is bounded in~${L^1(\Omega)}$ and then, up to a subsequence, weakly-${*}$ converges to some~${a \in \mathcal M(\Omega)}$. Moreover, the Radon measure~${a}$ belongs to~${L^\infty(\Omega)}$. Indeed, let~${\varphi\in\mathscr C_0(\Omega)}$ and denote again by~${\varphi}$ its extension to~${\R^3}$ by setting~${\varphi \equiv 0}$ on~${\R^3 \setminus \Omega}$. There exists a finite subset~${I_n}$ of~${\Z^3}$ such that 
\begin{equation*} 
\Omega \subset \bigcup_{k \in I_n} Q_n^k, 
\end{equation*}
where~${Q_n^k}$ is defined by~\eqref{decoup}. As~${\varphi}$ is a uniformly continuous function, we have 
\begin{equation} \label{lim1}
\int_\Omega \alpha_n \varphi \dx = \Sum_{k \in I_n} \int_{Q_n^k} \alpha_n \varphi \dx = \Sum_{k \in I_n} \varphi(\e_n k) \int_{Q_n^k} \alpha_n \dx  + {o} (1).
\end{equation}  
By~\eqref{ii}, we have 
\begin{equation} \label{lim2}
\Sum_{k \in I_n} |\varphi(\e_n k)| \int_{Q_n^k} \alpha_n \dx \leq  c \Sum_{k \in I_n} |Q_n^k| |\varphi(\e_n k)| = c \Vert \varphi \Vert_{L^1(\Omega)} + {o} (1).
\end{equation} 
The weak-$*$ convergence of~$\alpha_n$ to~$a$, combined with~\eqref{lim1} and~\eqref{lim2} yields
\begin{equation*} 
\left| \int_\Omega  \varphi \ \! a(\text{d}x)\right| \leq c \Vert \varphi \Vert_{L^1(\Omega)},
\end{equation*}
which implies that the measure~${a}$ is absolutely continuous with respect to the Lebesgue measure and~${a \in L^\infty(\Omega)}$. From~\eqref{convphas20} and the convergence of~${\alpha_n}$ to~${a}$, we have 
\begin{equation*}
\theta_n^{-1} \2 =  (\theta_n \alpha_{2,n})^{-1} (\alpha_n - \alpha_1 \1) \harpoon \theta:=\alpha_2^{-1} (a - \alpha_1) \in L^\infty(\Omega) \quad \text{weakly-${*}$ in } \mathcal M(\Omega),
\end{equation*} 
and then~${a = \alpha_1 + \theta \alpha_2}$.

On the other hand, by the Cauchy-Schwarz inequality combined with~\eqref{convlem},~\eqref{ii} and the convergence of~${\alpha_n}$ to~${\alpha_1 + \alpha_2 \, \theta}$, we have, for any~${\varphi \in \mathscr C_0(\Omega)}$, 
\begin{equation} \label{L2}
\begin{array}{l}
\vspace{0.4cm}\displaystyle\left|\int_\Omega \varphi \, \xi(\text{d}x) \right|^2 = \lim \limits_{n \to \infty} \left|\int_\Omega \xi_n \varphi \dx \right|^2 \leq \limsup \limits_{n \to \infty} \int_\Omega\alpha_n^{-1} \xi_n^2 \dx \displaystyle\int_\Omega\alpha_n \varphi^2 \dx\\
\displaystyle \leq c \, \int_\Omega(\alpha_1 + \theta \alpha_2) \varphi^2 \dx \leq c \Vert \alpha_1 + \theta \alpha_2 \Vert_\infty \ \Vert \varphi\Vert^2_{L^2(\Omega)},\end{array}
\end{equation} 
which implies that the limit measure~${\xi}$ of~${\xi_n}$ in~\eqref{convlem} is actually an element of~${L^2(\Omega)}$.

We now prove the convergence~\eqref{convdivcurl}. Let~${\varphi \in \mathscr D(\Omega)}$ and let~${I_n}$ be a finite subset of~${\Z^3}$ such that 
\begin{equation*}
\supp \varphi \subset \bigcup_{k \in I_n} \overline{Q}_n^{k} \subset \Omega,
\end{equation*} 
where~${\supp \varphi}$ is the support of~${\varphi}$. For any~${w \in H^1(\Omega)}$, define~${\overline{w}^{\e_n}}$ the piecewise constant function associated with the partition~${\left(Q_n^k \right)_{k \in I_n}}$ as follows: 
\begin{equation*}
\overline{w}^{\e_n} = \sum_{k \in I_n} \left( \fint_{Q_n^k} w\right) \mathds{1}_{Q_n^k}. 
\end{equation*} 
In order to study the convergence, in the sense of distributions, of~${(\xi_n v_n - \xi v)}$ to~${0}$, we rewrite it as 
\begin{equation} \label{div}
\xi_n v_n - \xi v  =  \underbrace{\xi_n \big(v_n - \overline{v_n}^{\e_n}\big)}_{:=p_n}  + \underbrace{\xi_n \big(\overline{v_n}^{\e_n} - \overline{v}^{\e_n}\big)}_{:=q_n}  +  \underbrace{\xi_n \overline{v}^{\e_n}\ - \xi v}_{:=r_n}
\end{equation} 
and estimate each term of the identity \eqref{div} separately.

\noindent \textit{Convergence of the term~${p_n}$ in~\eqref{div}}. Thanks to the Cauchy-Schwarz inequality, we have
\begin{equation} \label{bigineq}
\begin{array}{l}
 \displaystyle\left| \int_\Omega \xi_n \big(v_n - \overline{v_n}^{\e_n}\big) \varphi \dx \right|^2 \leq  \displaystyle\left(\sum_{k \in I_n} \left|\int_{Q_n^k} \xi_n \left(v_n - \fint_{Q_n^k} v_n \right) \varphi \dx \right|\right)^2 \vspace{0.2cm} \\
\leq  \displaystyle \Vert \varphi \Vert^2_\infty \left(\sum_{k \in I_n}\sqrt{\int_{Q_n^k} \alpha_n^{-1} |\xi_n|^2 \dx }  \vspace{0.3cm}\sqrt{\int_{Q_n^k}   \alpha_n \left| v_n - \fint_{Q_n^k} v_n \right|^2  \dx } \right)^2\\
\leq  \displaystyle \Vert \varphi \Vert^2_\infty \sum_{k \in I_n} \int_{Q_n^k} \alpha_n^{-1} |\xi_n|^2 \dx \vspace{0.3cm}\sum_{k \in I_n} \int_{Q_n^k}   \alpha_n \left| v_n - \fint_{Q_n^k} v_n \right|^2  \dx \\
\leq  \displaystyle c_n \Vert \varphi \Vert^2_\infty {\int_{\Omega} \alpha_n^{-1} |\xi_n|^2 \dx }  \vspace{0.3cm}\  { \int_{\Omega}   \alpha_n \left|  \n v_n \right|^2  \dx }, \end{array} 
\end{equation} 
where the last inequality is a consequence of~\eqref{WPW}. Finally, the inequality~\eqref{bigineq} combined with~\eqref{divrot1} and the convergence~\eqref{condWPW} yield \vspace{-0.2cm}
\begin{equation} \label{convlem2}
\displaystyle\left| \int_\Omega \xi_n \big(v_n - \overline{v_n}^{\e_n}\big) \varphi \dx \right| \leq c \sqrt{c_n} \xrightarrow[n \to \infty]{} 0. 
\end{equation}

\noindent \textit{Convergence of the term~${q_n}$ in~\eqref{div}.} By the Cauchy-Schwarz inequality and~\eqref{ii}, we have 
\begin{equation*} 
\begin{array}{l}
\displaystyle\left| \int_\Omega \xi_n \big(\overline{v_n}^{\e_n} - \overline{v}^{\e_n}\big) \varphi \dx \right|^2 \leq \vspace{0.4cm} \displaystyle\left(\sum_{k \in I_n} \left|\int_{Q_n^k} \xi_n \ \varphi \dx \right| \left| \fint_{Q_n^k} (v_n - v) \dx \right|\right)^2 \\
 \leq  \vspace{0.4cm}\displaystyle  \Vert \varphi \Vert^2_\infty \left(\sum_{k \in I_n} |Q_n^k|^{-1} \sqrt{\int_{Q_n^k} \alpha_n \dx} \sqrt{\int_{Q_n^k} \alpha_n^{-1} |\xi_n|^2 \dx } \vspace{0.4cm} \int_{Q_n^k}  |v_n - v| \dx \right)^2 \\
 \leq  \displaystyle \vspace{0.4cm} \Vert \varphi \Vert^2_\infty \left(\sum_{k \in I_n} \sqrt{\fint_{Q_n^k} \alpha_n \dx} \sqrt{\int_{Q_n^k} \alpha_n^{-1} |\xi_n|^2 \dx } \vspace{0.4cm}\sqrt{\int_{Q_n^k}  (v_n - v)^2 \dx}\right)^2 \\
  \leq \displaystyle \vspace{0.1cm}c \Vert \varphi \Vert^2_\infty {\sum_{k \in I_n} \int_{Q_n^k} \alpha_n^{-1} |\xi_n|^2 \dx } {\sum_{k \in I_n} \int_{Q_n^k}  (v_n - v)^2 \dx} \\
  \leq \displaystyle \vspace{0.2cm}c {\int_{\Omega} \alpha_n^{-1} |\xi_n|^2 \dx} \ { \int_{\Omega} (v_n - v)^2 \dx } 
\end{array}
\end{equation*}
which yields, by~\eqref{divrot1},
\begin{equation} \label{ineqlong}
\displaystyle\left| \int_\Omega \xi_n \big(\overline{v_n}^{\e_n} - \overline{v}^{\e_n}\big) \varphi \dx \right|^2   \leq \displaystyle c \, \Vert v_n - v \Vert^2_{L^2(\Omega)}. 
\end{equation} 
Since~${v_n}$ converges weakly to~${v}$ in~${H_0^1(\Omega)}$, by Rellich's theorem, up to a subsequence,~${v_n}$ converges strongly to~${v}$ in~${L^2(\Omega)}$. Hence,~\eqref{ineqlong} implies that 
\begin{equation} \label{convlem1}
\int_\Omega \xi_n \big(\overline{v_n}^{\e_n} - \overline{v}^{\e_n}\big) \varphi \dx  \xrightarrow[n \to \infty]{} 0.
\end{equation}

\noindent \textit{Convergence of the term~${r_n}$ in~\eqref{div}.} Consider, for any~${\delta >0}$, an approximation~${\psi_\delta \in \mathscr C_c({\Omega})}$ of~${v}$ for the~${L^2(\Omega)}$ norm, \textit{i.e.}, 
\begin{equation} \label{approx}
\Vert v - \psi_\delta\Vert_{L^2(\Omega)} = o(\delta).
\end{equation} 
The term~${r_n}$ in~\eqref{div} writes 
\begin{equation} \label{decomp''} 
\xi_n \overline{v}^{\e_n} - \xi v = \xi_n \big(\overline{v}^{\e_n} - \overline{\psi_\delta}^{\e_n} \big) + \xi_n \big( \overline{\psi_\delta}^{\e_n} - \psi_\delta \big) + (\xi_n - \xi) \psi_\delta + \xi (\psi_\delta - v).
\end{equation} 
On the one hand, since~${\overline{\psi_\delta}^{\e_n}}$ converges uniformly, as~${n}$ goes to infinity, to~${\psi_\delta \in \mathscr C_c ({\Omega})}$, the convergence~\eqref{convlem} of~${\xi_n}$  to~${\xi}$ implies that the second term and the third term in the right hand side of the equality~\eqref{decomp''} converge to~${0}$ in~${\mathscr D'(\Omega)}$. Moreover, by the Cauchy-Schwarz inequality and the fact that~${\xi \in L^2(\Omega)}$, we have 
\begin{equation*}
\left|\int_\Omega \xi (\psi_\delta - v) \varphi \dx \right| \leq \Vert \varphi \Vert_\infty \Vert \xi \Vert_{L^2(\Omega)} \Vert \psi_\delta - v \Vert_{L^2(\Omega)}.
\end{equation*} On the other hand, following~\eqref{ineqlong}, we have 
\begin{equation} \label{conv32} 
\displaystyle\left| \int_\Omega \xi_n \big(\overline{v}^{\e_n} - \overline{\psi_\delta}^{\e_n}\big) \varphi \dx \right| \leq c \, \Vert \varphi \Vert_\infty \Vert v - \psi_\delta \Vert_{L^2(\Omega)}.
\end{equation} 
Hence, by~\eqref{approx}-\eqref{conv32}, we have
\begin{equation} \label{convfin'} \limsup \limits_{n \to \infty} \left| \int_\Omega \big(\xi_n \overline{v_n}^{\e_n} - \xi v \big) \varphi \dx \right| \leq c \, \Vert \psi_\delta - v \Vert_{L^2(\Omega)} = o(\delta), 
\end{equation} 
for arbitrary~${\delta >0}$.

Finally, putting together~\eqref{div},~\eqref{convlem2},~\eqref{convlem1} and~\eqref{convfin'}, we obtain that 
\begin{equation*} 
\limsup \limits_{n \to \infty} \left| \int_\Omega (\xi_n v_n - \xi v) \varphi \dx \right| = o(\delta),
\end{equation*} 
which concludes the proof of Lemma~\ref{divcurl}. \qed

\medskip

\noindent In the sequel we apply Lemma~\ref{divcurl} to sequences~${\xi_n}$ of vector-valued functions in~${L^1(\Omega)^2}$ or~${L^1(\Omega)^3}$.

\bigskip

\noindent \textbf{{Proof of Theorem~\ref{thgen}.}} Thanks to the equi-coerciveness~${\sigma_n \geq \alpha_1 I_3}$, the solution~${u_n}$ of the problem~\eqref{Dir} satisfies the convergence, up to a subsequence,
\begin{equation} \label{convu}
u_n \harpoon u \quad \text{weakly in } H_0^1(\Omega),
\end{equation} 
for some~${u}$ in~${H_0^1(\Omega)}$. Moreover, putting~${u_n}$ as a test function in the equation~\eqref{Dir}, we obtain that 
\begin{equation} \label{ineqener}
\int_\Omega \alpha_n |\nabla u_n|^2 \dx =  \int_\Omega \sigma_n \nabla u_n \cdot \nabla u_n \dx = \langle f, u_n\rangle_{H^{-1}(\Omega), H_0^1(\Omega)} \leq c.
\end{equation} 
Since~${\alpha_{2,n}^{-1} \sigma_{2,n} = I_3 + \alpha_{2,n}^{-1} \, \beta_{2,n} \mathscr E(h)}$, by~\eqref{convphas20} the sequence~${|\alpha_{2,n}^{-1} \sigma_{2,n}|}$ is bounded. Then, as the sequence~${\alpha_n}$ is bounded in~${L^1(\Omega)}$, the Cauchy-Schwarz inequality and~\eqref{ineqener} give 
\begin{equation*} 
\left(\int_\Omega |\sigma_n \nabla u_n|\dx\right)^2 \leq c {\int_\Omega \alpha_n \dx} {\int_\Omega \alpha_n |\nabla u_n|^2 \dx } \leq c.
\end{equation*} 
Hence, we have the convergence of the current~${\sigma_n \n u_n}$, up to a subsequence, 
\begin{equation} \label{limcour}
\sigma_n \nabla u_n \harpoon \xi_0 \quad \text{weakly-${*}$ in } \mathcal M(\Omega)^3, 
\end{equation} 
for some~${\xi_0 \in \mathcal M(\Omega)^3}$. Moreover, by the boundedness of~${|\alpha_{2,n}^{-1} \sigma_{2,n}|}$ and~\eqref{ineqener}, we have 
\begin{equation*}
\int_\Omega \alpha_n^{-1} |\sigma_n \nabla u_n|^2 \dx \leq c \int_\Omega \alpha_n |\nabla u_n|^2 \dx \leq c.
\end{equation*}
Then, by Lemma~\ref{divcurl} applied to~${\xi_n:= \sigma_n \n u_n}$, the measure~${\xi_0}$ is actually an element of~${L^2(\Omega)^3}$. 

\bigskip

The rest of the proof, which is divided into three steps, is devoted to the determination of the form of the limit current~$\xi_0$. To that end, we use a method in the spirit of~\mbox{H-convergence} of Murat-Tartar which is based on the cylindrical nature of the microstructure and the compactness result of Lemma~\ref{divcurl} for sequences only bounded in~${L^2(\Omega ; \sigma_n^{-1/2} \! \! \dx)}$. In the first two steps, we compute the components~$\xi_0 \cdot e_1$ and~$\xi_0 \cdot e_2$  by combining Lemma~\ref{divcurl} with a corrector function associated with the transversal conductivity~${\w \sigma_n}$, the existence of which is ensured by the two-dimensional results of~\mbox{\cite{NHM,BP1}}. Since the corrector function considered in the previous steps is independent of the variable~${x_3}$, the component~${\xi_0 \cdot e_3}$ needs a different approach. This is the object of the last step.

\bigskip

\noindent \textit{First step:} Building a corrector. Thanks to Remark~\ref{remexpl'}, up to a subsequence,~${\w \sigma_n}$~${H(\mathcal M(\w \Omega)^2)}$-converges to some coercive matrix-valued function~${\w \sigma_*}$. Then, the sequence~${\w \sigma_n^\T}$ \mbox{${H(\mathcal M(\w \Omega)^2)}$-converges} to~${\w \sigma_*^\T}$~(see Theorem 2.1 of~\cite{NHM}). Let~${\lambda \in \R^3}$ with~${\lambda \perp e_3}$. For~${\w \lambda = (\lambda_1,\lambda_2)^\T \in \R^2}$, let~${v_n^{\w \lambda}}$ be the solution of 
\begin{equation} \label{quas:Rn}
\left\{\begin{array}{r c l l}
\divg \big(\w \sigma_n^\T \w \nabla v_n^{\w \lambda} \big) & = & \mathrm{div} \big(\w \sigma_*^\T \w \lambda\big)\vspace{0.2cm} &\text{in } \w \Omega\\
v_n^{\w \lambda } & = & \w \lambda \cdot \w x & \text{on } \partial \w \Omega.\\
\end{array} \right.
\end{equation} 
By Definition~\ref{HM}, we have the convergences 
\begin{equation} \label{convvn}
\left\{ \begin{array}{ r l l}
v_n^{ \w \lambda} & \harpoon \w \lambda \cdot \w x = \lambda \cdot x \vspace{0.2cm} & \text{weakly in } H^1(\w \Omega), \\
\w \sigma_n^\T \w \nabla v_n^{\w \lambda} & \harpoon \w \sigma_*^\T \w \lambda &\text{weakly-${*}$ in }\mathcal M (\w \Omega)^2.\\
\end{array} \right.\end{equation} 
Setting, for~${x \in \Omega}$,~${v_n^\lambda(x_1,x_2,x_3) = v_n^{ \w \lambda}(x_1,x_2)}$, we have the convergences 
\begin{equation} \label{convcorr}
\left\{\begin{array}{ r l l}
v_n^{ \lambda} & \harpoon \lambda \cdot x \vspace{0.2cm} & \text{weakly in } H^1(\Omega), \\
\w \sigma_n^\T \w \nabla v_n^{ \lambda} & \harpoon \w \sigma_*^\T \w \lambda &\text{weakly-${*}$ in }\mathcal M (\Omega)^2,\\
\end{array} \right.
\end{equation} 
and the energy inequality, as in~\eqref{ineqener}, 
\begin{equation} \label{ineqenercorr}
\int_\Omega \alpha_n \big|\nabla v_n^{\lambda}\big|^2 \dx = \int_{\w \Omega} \w \sigma_n \w \nabla v_n^{\w \lambda} \cdot \w \nabla v_n^{\w \lambda} \ \text{d}\w x  \leq c.
\end{equation}

Let~${\varphi \in \mathscr D(\Omega)}$. By~\eqref{limcour} and since~${v_n^\lambda}$ converges weakly to~${\lambda \cdot x}$ in~${H^1(\Omega)}$, putting~${v_n^{ \lambda} \varphi}$ as a test function in~\eqref{Dir} yields 
\begin{equation} \label{quas:convbase}
\int_\Omega \sigma_n \nabla u_n \cdot \nabla \big( v_n^{ \lambda} \varphi\big) \dx \xrightarrow[n \to \infty]{} \langle f, \varphi \ \lambda \cdot x \rangle_{H^{-1}(\Omega),H_0^1(\Omega)} = \int_\Omega \xi_0 \cdot \n(\varphi \lambda \cdot x) \dx.
\end{equation} 
Since~${\sigma_n}$ and~${v_n^{\lambda}}$ do not depend on the variable~${x_3}$, we have the identity 
\begin{equation} \label{quas:decomp2} 
\sigma_n \nabla u_n \cdot \nabla v_n^{\lambda} = \w \sigma_n \w \nabla u_n \cdot \w \nabla v_n^{\lambda} - \partial_3 \big( \beta_n \w \nabla v_n^{\lambda} \cdot J \w h \ u_n\big). 
\end{equation} 
Then, by~\eqref{quas:decomp2},
an integration by parts gives 
\begin{eqnarray}
\displaystyle \int_\Omega \sigma_n \nabla u_n \cdot \nabla \big( v_n^{\lambda} \varphi\big) \dx & = & \displaystyle\int_\Omega \sigma_n \nabla u_n \cdot \nabla \varphi \ v_n^{\lambda} \dx \label{decompcomp1} \\ 
&\displaystyle + & \displaystyle \int_\Omega \beta_n  \w \nabla v_n^{\lambda} \cdot J \w h \ u_n \ \frac{\partial \varphi}{\partial x_3} \dx \label{decompcomp2}\\&\displaystyle +& \displaystyle \int_\Omega \w \sigma_n^\T \w \nabla v_n^{\lambda} \cdot \w \nabla u_n \ \varphi \dx. \label{decompcomp3}
\end{eqnarray}

\bigskip

\noindent \textit{Step 2:} Estimates of the terms in~\eqref{decompcomp1}-\eqref{decompcomp3}. The convergence of these terms are consequences of Lemma~\ref{divcurl} and the generalized two-dimensional div-curl lemma in a high-contrast context of~\cite{BrCa1}.

\medskip

\noindent \textit{Convergence of the term on the right hand side of~\eqref{decompcomp1}.} On the one hand, by the boundedness of~${\alpha_n^{-1} \sigma_n}$ and~\eqref{ineqener}, we have the inequality 
\begin{equation*} 
\int_\Omega \alpha_n^{-1}|\sigma_n \nabla u_n \cdot \nabla \varphi |^2 \dx \leq c \, \Vert \n \varphi \Vert_\infty^2 \int_\Omega \alpha_n |\nabla u_n|^2 \leq c.
\end{equation*} 
On the other hand, the convergence~\eqref{limcour}, the inequality~\eqref{ineqenercorr}, and the convergence~\eqref{convcorr} of~${v_n^\lambda}$ to~${\lambda \cdot x}$, show that the sequences~${\xi_n:=\sigma_n \nabla u_n}$ and~${v_n := v_n^{\lambda}}$ satisfy the assumptions of Lemma~\ref{divcurl}. Hence, we obtain 
\begin{equation} \label{decompcomp1'}
\displaystyle\int_\Omega \sigma_n \nabla u_n \cdot \nabla \varphi \ v_n^{ \lambda} \dx \xrightarrow[n \to \infty]{} \int_\Omega \xi_0 \cdot \nabla \varphi \ \lambda \cdot x \dx.
\end{equation}
 
\medskip

\noindent \textit{Convergence of the term in~\eqref{decompcomp2}.} We first compute the limit of~${\beta_n \w \n v_n^{\lambda}}$ in the sense of Radon measures. We have the identity 
\begin{equation} \label{identity}
\begin{array}{r c l} \beta_n \w \nabla v_n^{\lambda} = \beta_n \left(  \frac{\partial v_n^\lambda}{\partial x_1},\frac{\partial v_n^\lambda}{\partial x_2}\right)^\T & = & \1 \beta_1 \w \nabla v_n^{\lambda} + \beta_{2,n} \w \sigma_{2,n}^{-T}\left[ \2 \w \sigma_{2,n}^\T \w \n v_n^\lambda \right] \\*[0.9em]
& = & \1 \beta_1 \w \nabla v_n^{\lambda} + \beta_{2,n} \w \sigma_{2,n}^{-T}\left[ \w \sigma_n^\T \w \nabla v_n^{\lambda} - \1 \w \sigma_1^\T \w \nabla v_n^{\lambda} \right], \end{array}
\end{equation} 
where 
\begin{equation} \label{identity2}
\w \sigma_1 := \begin{pmatrix} \alpha_1 & - \beta_1 \, h_3 \\*[0.3em]
\beta_1 \, h_3 & \alpha_1 \end{pmatrix} \quad \text{and} \quad  \w \sigma_{2,n} := \begin{pmatrix} \alpha_{2,n} & - \beta_{2,n} \, h_3 \\*[0.3em]
\beta_{2,n} \, h_3 & \alpha_{2,n} \end{pmatrix}.
\end{equation} 
By~\eqref{convphas20}, we have 
\begin{equation} \label{convdir}
\beta_{2,n} \w \sigma_{2,n}^{-T} \xrightarrow[n \to \infty]{} \beta_{2} \w \sigma_2^{-\T}, \quad \text{where} \quad \w \sigma_{2} := \begin{pmatrix} \alpha_{2} & - \beta_{2} \, h_3 \\*[0.3em]
\beta_{2} \, h_3 & \alpha_{2} \end{pmatrix}.
\end{equation}
Combining this convergence with the ones in~\eqref{convvn}, we obtain that 
\begin{equation} \label{convbid} 
\beta_n \w \nabla v_n^{\lambda} \harpoon \beta_1 \w \lambda + \beta_2 \w \sigma_{2}^{-\T} \left[ \w \sigma_*^\T - \w \sigma_1^\T\right] \w \lambda \quad \text{weakly-${*}$ in } \mathcal M(\Omega)^2. 
\end{equation} 
By the boundedness of~${\alpha_n^{-1} \sigma_n}$,~\eqref{convphas20} and~\eqref{ineqenercorr}, we have 
\begin{equation*} 
\int_\Omega \alpha_n^{-1}|\beta_n \w \nabla v_n^{\lambda}|^2 \dx \leq c \int_\Omega \alpha_n |\nabla v_n^{ \lambda}|^2 \leq c.
\end{equation*} 
This inequality together with~\eqref{ineqener},~\eqref{convbid} and the weak convergence~\eqref{convu} of~${u_n}$ to~${u}$ in~${H_0^1(\Omega)}$ show that the sequences~${{\xi_n:=\beta_n \w \nabla v_n^{\lambda}}}$ and~${v_n := u_n}$ satisfy the assumptions of Lemma~\ref{divcurl}. Then, 
\begin{equation} \label{decompcomp2'}
\int_\Omega \beta_n \w \nabla v_n^{\lambda} \cdot J\w h \ u_n \ \frac{\partial \varphi}{\partial x_3} \dx \xrightarrow[n \to \infty]{} \int_\Omega  \left(\beta_1 I_2 + \beta_2 \left[ \w \sigma_* - \w \sigma_1 \right] \w \sigma_{2}^{-1} \right) J \w h \cdot \w \lambda \ u \ \frac{\partial \varphi}{\partial x_3} \dx.
\end{equation}

\noindent \textit{Convergence of the term in~\eqref{decompcomp3}.} Integrating by parts in~\eqref{decompcomp3}, we obtain that 
\begin{equation} \label{bid1}
\int_\Omega \w \sigma_n^\T \w \n v_n^{\lambda } \cdot \w \n u_n \ \varphi \dx = - \int_\Omega \w \sigma_n^\T \w \n v_n^{\lambda } \cdot \w \n \varphi \ u_n \dx + \int_\Omega \w \sigma_n^\T \w \n v_n^{\lambda } \cdot  \w \nabla \left(u_n \varphi \right) \dx .
\end{equation}
On the one hand, the boundedness of~${\alpha_n^{-1} \sigma_n}$ and~\eqref{ineqenercorr} yields 
\begin{equation*}
\int_\Omega \alpha_n^{-1} \big|\w \sigma_n^\T \w \n v_n^{\lambda } \big|^2 \leq c \int_\Omega \alpha_n \big|\w \n v_n^{\lambda } \big|^2 \leq c.
\end{equation*} 
By the second convergence of~\eqref{convcorr}, the weak convergence~\eqref{convu} of~${u_n}$ to~${u}$ in~${H_0^1(\Omega)}$ and~\eqref{ineqener}, the sequences~${\xi_n := \w \sigma_n^\T \w \n v_n^{\lambda }}$ and~${v_n:=u_n}$ satisfy the assumptions of Lemma~\ref{divcurl}. Hence, 
\begin{equation} \label{bid1'}
\int_\Omega \w \sigma_n^\T \w \n v_n^{\lambda } \cdot \w \n \varphi \ u_n \dx \xrightarrow[n \to \infty]{} \int_\Omega \w \sigma_*^\T {\w \lambda } \cdot \w \n \varphi \ u \dx.
\end{equation} 
On the other hand, since~${\w \sigma_n^\T \w \n v_n^{\w \lambda }}$ does not depend on the variable~${x_3}$, the second term on the right hand side of~\eqref{bid1} can be rewritten under the form 
\begin{equation} \label{bid''}
\int_\Omega \w \sigma_n^\T \w \n v_n^{\lambda } \cdot  \w \nabla \left(u_n \varphi \right) \dx = \int_{\w \Omega} \w \sigma_n^\T \w \n v_n^{\w \lambda }(x') \cdot  \w \nabla \left[ \int_0^1\left(u_n \varphi \right)(x',x_3) \ \text{d}x_3 \right] \ \text{d}x'
\end{equation} 
where~${x' = (x_1,x_2)}$. 

In order to study the asymptotic behavior of~\eqref{bid''}, we apply a two-dimensional div-curl lemma of~\cite{BP1} 
which is an extension to the non symmetric case of~\cite{BrCa1}. Set, for any~${x' \in \Omega}$, \begin{equation} \label{defeta} 
\eta_n(x') := \w \sigma_n^\T \w \n v_n^{\w \lambda }(x') \quad \text{and} \quad v_n(x') : = \int_0^1\left(u_n \varphi \right)(x',x_3) \ \text{d}x_3. 
\end{equation} 
Due to the convergences~\eqref{convvn} and~\eqref{convu}, we have 
\begin{equation} \label{convbid'''}
\left\{ \begin{array}{r c l l} \eta_n  & \harpoon & \w \sigma_*^\T \w \lambda & \text{weakly-${*}$ in } \mathcal M(\w \Omega)^2\\*[0.3em] v_n(x') & \harpoon & \displaystyle v(x'):=\int_0^1\left(u \varphi \right)(x',x_3) \ \text{d}x_3 & \text{weakly in } H^1(\w \Omega). \end{array}\right.
\end{equation}
The convergences~\eqref{convphas20} and~\eqref{tt'}, the definition~\eqref{quas:Rn} of the corrector~${v_n^{\w \lambda}}$,~\eqref{convvn}-\eqref{ineqenercorr} and~\eqref{convbid'''} imply that the sequences~${\eta_n}$ and~${v_n}$ defined in~\eqref{defeta} satisfy the assumptions of Lemma 2.1 in~\cite{BP1}. Then, 
\begin{equation} \label{divrotbid}
\w \sigma_n^\T \w \n v_n^{\w \lambda }(x') \cdot  \w \nabla \left[ \int_0^1\left(u_n \varphi \right)(x',x_3) \ \text{d}x_3 \right] \harpoon \w \sigma_*^\T {\w \lambda} \cdot  \w \nabla \left[ \int_0^1\left(u \varphi \right)(x',x_3) \ \text{d}x_3 \right] \quad \text{in } \mathscr D'(\w \Omega).
\end{equation} 
Let~${\psi \in \mathscr D(\w \Omega)}$ such that~${\psi \equiv1}$ on the projection of the support of~${\varphi}$ on the~${(x_1,x_2)}$-plane. Taking~${\psi}$ as a test function in~\eqref{divrotbid}, we obtain 
\begin{equation*} 
\displaystyle \int_{\w \Omega} \w \sigma_n^\T \w \n v_n^{\w \lambda }(x') \cdot  \w \nabla \left[ \int_0^1\left(u_n \varphi \right)(x',x_3) \ \text{d}x_3 \right] \ \text{d}x' \xrightarrow[n\to \infty]{} \int_\Omega \sigma_*^\T {\w \lambda} \cdot \w \n  \left( u \varphi \right).
\end{equation*} 
Finally, this convergence combined with~\eqref{bid1},~\eqref{bid1'} and~\eqref{bid''} gives 
\begin{equation} \label{conv3}
\int_\Omega \w \sigma_n^\T \w \n v_n^{ \lambda } \cdot \w \n u_n \ \varphi \dx \xrightarrow[n \to \infty]{} \int_\Omega  \w \sigma_*^\T {\w \lambda } \cdot \w \nabla u \ \varphi \dx .  
\end{equation} 
Putting together~\eqref{decompcomp1'},~\eqref{decompcomp2'} and ~\eqref{conv3} with the equality~\eqref{decompcomp1}-\eqref{decompcomp3}, we obtain that 
\begin{equation*}
\begin{array}{r l}
\hspace{-0.5cm}\vspace{0.4cm}\displaystyle \int_\Omega \sigma_n \nabla u_n \cdot \nabla \big( v_n^{ \lambda} \varphi\big) \dx & \xrightarrow[n \to \infty]{} \displaystyle\int_\Omega \xi_0 \cdot \nabla \varphi \ \lambda \cdot x \dx + \int_\Omega \w \sigma_* \w \n u \cdot {\w \lambda} \ \varphi \dx \\ 
&\displaystyle \hspace{1.5cm} + \int_\Omega \left( \beta_1 I_2 + \beta_2 \left[ \w \sigma_* - \w \sigma_1 \right] \w \sigma_{2}^{-1} \right) J \w h\cdot \w \lambda \ u \ \frac{\partial \varphi}{\partial x_3} \dx.
\end{array} 
\end{equation*} 
Since~${\w \sigma_*}$ depends only on the variable~${(x_1,x_2)}$, this convergence, an integration by parts and~\eqref{quas:convbase} give 
\begin{equation} \label{t't} 
\int_\Omega \xi_0 \cdot \lambda \ \varphi \dx = \int_\Omega \left[\w \sigma_* \w \n u - \frac{\partial u}{\partial x_3} \left( \beta_1 I_2 +  \beta_2 \left[ \w \sigma_* - \w \sigma_1 \right] \w \sigma_{2}^{-1} \right) J \w h \right] \cdot \w \lambda \  \varphi \dx. 
\end{equation} 
Finally, since the equation~\eqref{t't} holds for any~${\varphi \in \mathscr D(\Omega)}$ and any~${\lambda \perp e_3}$, we obtain the first two components of~${\xi_0}$ 
\begin{equation} \label{limcour'} 
\w \xi_0 = \w \sigma_* \w \n u - \displaystyle \frac{\partial u}{\partial x_3} \left( \beta_1 I_2 +  \beta_2 \left[ \w \sigma_* - \w \sigma_1 \right] \w \sigma_{2}^{-1} \right) J \w h.
\end{equation}
 
\noindent \textit{Step 3:} Computation of~${\xi_0 \cdot e_3}$. By~\eqref{limcour}, we have the convergence 
\begin{equation} \label{limdeb}
\alpha_n \, \frac{\partial u_n}{\partial x_3} + \beta_n \w \n u_n \cdot J \w h = \sigma_n \n u_n \cdot e_3\harpoon \xi_0 \cdot e_3 \quad \text{weakly-${*}$ in } \mathcal M(\Omega).
\end{equation} 

We first study the asymptotic behaviour of~${\alpha_n \, \partial_3 u_n}$~(which also gives the limit of~${\beta_n \, \partial_3 u_n}$ due to the fact that, by virtue of \eqref{convphas20},~${\alpha_{2,n}}$ and~${\beta_{2,n}}$ are of the same order). On the one hand, since~${\theta_n = |\Omega|^{-1} |\Omega_n|}$, by the convergence~\eqref{convphas20}, we have 
\begin{equation*} 
\int_\Omega \alpha_n^{-1} \ | \theta_n^{-1}\2|^2 \dx =  \frac{\theta_n^{-1}|\Omega_n|}{\theta_n \alpha_{2,n}} = \frac{|\Omega|}{\theta_n \alpha_{2,n}} \leq  c.
\end{equation*} 
On the other hand, by~\eqref{ineqener}, the weak convergence~\eqref{convu} of~${u_n}$ to~${u}$  in~${H_0^1(\Omega)}$ and~\eqref{tt}, the sequences~${\xi_n:= \theta_n^{-1}\2}$ and~${v_n:=u_n}$ satisfy, once again, the assumptions of Lemma~\ref{divcurl}. Hence, 
\begin{equation*}
\theta_n^{-1}\2 u_n \harpoon \theta u \quad \text{in } \mathscr D'(\Omega).
\end{equation*}
Moreover, since~${\2}$ does not depend on the variable~${x_3}$, we have 
\begin{equation} \label{limfib}
\theta_n^{-1}\2 \frac{\partial u_n}{\partial x_3} \harpoon \theta \frac{\partial u}{\partial x_3} \quad \text{in } \mathscr D'(\Omega).
\end{equation} 
Finally, thanks to~\eqref{limfib} and~\eqref{convphas20}, we obtain the convergences, in the sense of Radon measures, 
\begin{equation} \label{limfincur}
\left\{\begin{array}{r c l} \displaystyle\alpha_n \ \frac{\partial u_n}{\partial x_3} = \vspace{0.4cm}\alpha_1 \1 \ \frac{\partial u_n}{\partial x_3} + (\theta_n \alpha_{2,n}) \theta_n^{-1} \2 \ \frac{\partial u_n}{\partial x_3} &\harpoon & \displaystyle  (\alpha_1 + \theta \alpha_2)\frac{\partial u}{\partial x_3},\\
\displaystyle\beta_n \ \frac{\partial u_n}{\partial x_3} = \beta_1 \1 \ \frac{\partial u_n}{\partial x_3} + (\theta_n \beta_{2,n}) \theta_n^{-1} \2 \ \frac{\partial u_n}{\partial x_3} &\harpoon & \displaystyle (\beta_1 + \theta \beta_2)\frac{\partial u}{\partial x_3}.\end{array}\right.
\end{equation} 
Now, in order to obtain the limit of the term~${\beta_n \w \n u_n}$ in~\eqref{limdeb}, which similarly to~\mbox{\eqref{identity}-\eqref{identity2}}, writes
\begin{equation} \label{tttt}
\beta_n \w \n u_n = \beta_1 \1 \w \n u_n + \beta_{2,n} \w \sigma_{2,n}^{-1} \left[ \w \sigma_n \w \n u_n - \1 \w \sigma_1 \w\n u_n\right],
\end{equation} 
it remains to estimate~${\w \sigma_n \w \n u_n}$. Since~${\xi_0}$ is the limit of the current~${\sigma_n \n u_n}$~\eqref{limcour} and since 
\begin{equation*} 
\forall \lambda \perp e_3, \quad \sigma_n \n u_n \cdot \lambda = \w \sigma_n \w \n u_n \cdot \w \lambda - \beta_n \ \displaystyle \frac{\partial u_n}{\partial x_3} J \w h \cdot \w \lambda,
\end{equation*}
the equality~\eqref{limcour'} gives 
\begin{equation*} 
\w \sigma_n \w \n u_n - \beta_n \ \displaystyle \frac{\partial u_n}{\partial x_3} J \w h\harpoon \w \sigma_* \w \n u - \frac{\partial u}{\partial x_3} \left( \beta_1  I_2 + \beta_2 \left[ \w \sigma_* - \w \sigma_1 \right] \w \sigma_{2}^{-1} \right) J \w h \quad \text{weakly-${*}$ in } \mathcal M(\Omega)^2.
\end{equation*} 
Then, combining this convergence with~\eqref{limfincur}, we have 
\begin{equation} \label{lim''} 
\w \sigma_n \w \n u_n \harpoon  \w \sigma_* \w \n u + \beta_2 \, \frac{\partial u}{\partial x_3} \left( \left[\w \sigma_1 + \theta \w \sigma_2 - \w \sigma_* \right] \w \sigma_{2}^{-1} \right) J \w h \quad \text{weakly-${*}$ in } \mathcal M(\Omega)^2. 
\end{equation} 
Finally, passing to the limit in~\eqref{tttt}, taking into account~\eqref{lim''} and~\eqref{convdir}, we obtain the convergence, in the sense of Radon measures, 
\begin{equation} \label{fineq}
\beta_n \w \nabla u_n \harpoon \left[ \beta_1 I_2 + \beta_2 \, \w \sigma_2^{-1} (\w \sigma_* - \w\sigma_1)  \right] \w \n u + \beta_2^2 \ \frac{\partial u}{\partial x_3} \w \sigma_2^{-1} \big(\w \sigma_1 + \theta \w \sigma_2  - \w \sigma_*\big) \w \sigma_2^{-1} J \w h. 
\end{equation} 
Putting together~\eqref{limdeb},~\eqref{limfincur} and~\eqref{fineq} yields
\begin{equation} \label{limcour''}
\begin{array}{r l}
\xi_0 \cdot e_3 = & \left[ \beta_1 I_2 + \beta_2 \w \sigma_2^{-1} (\w \sigma_* - \w\sigma_1)  \right]^\T J \w h \cdot \w \n u \\*[0.2cm]
 & + \left[(\alpha_1 + \alpha_2 \theta) +  \beta_2^2 \w \sigma_2^{-1} \big(\w \sigma_1 + \theta \w \sigma_2  - \w \sigma_*\big) \w \sigma_2^{-1} J \w h \cdot J \w h \right] \displaystyle \frac{\partial u}{\partial x_3}.
\end{array}
\end{equation}
Finally, since the current~${\sigma_n \n u_n}$ weakly-${*}$ converges to~${\xi_0}$ in~\eqref{limcour}, we have the limit equation 
\begin{equation*} 
{- \divg(\xi_0) = f}, 
\end{equation*} 
where, by~\eqref{limcour'} and~\eqref{limcour''}, 
\begin{equation*} 
\xi_0  =  ( \w \xi_0, \xi_0 \cdot e_3)^\T = \sigma_*(h) \n u 
\end{equation*} 
which yields to the expression~\eqref{propgen'}-\eqref{eqgen} of~${\sigma_*(h)}$. Theorem~\ref{thgen} is proved. 
\qed

\markboth{Chapitre~\ref*{chap2} - Homogenization of high-contrast and non symmetric conductivities for non periodic columnar structures}{\textit{Section~\ref*{examples} - Two examples}}


\section{Two examples} \label{examples}

\markboth{Chapitre~\ref*{chap2} - Homogenization of high-contrast and non symmetric conductivities for non periodic columnar structures}{\textit{Section~\ref*{examples} - Two examples}}

In this section we present two examples where the perturbation formulas for the effective conductivities of non periodic high-contrast columnar composites are fully explicitly computed. 


\subsection{Circular fibres with variable radius}

\begin{figure}[H]
\centering
\begin{tikzpicture}[scale=1.7]
\clip (-3.5,-1.1) rectangle (3,2);

\foreach \col in {-0.5,-0.25,0,0.25,0.5,0.75,1,1.25,1.5,1.75}{
\draw (-1.5,\col) -- (3,\col);}

\foreach \col in {-1.5,-1.25,-1,-0.75,-0.5,-0.25,0,0.25,0.5,0.75,1,1.25,1.5,1.75,2,2.25,2.5,2.75,3}{
\draw (\col,-1) -- (\col,1.95);}
 
 \foreach \col in {-1.5+0.125,-1.25+0.125,-1+0.125,-0.75+0.125,-0.5+0.125,-0.25+0.125,0+0.125,0.25+0.125,0.5+0.125,0.75+0.125,1+0.125,1.25+0.125,1.5+0.125,1.75+0.125,2+0.125}{\foreach \hi in {-1+0.125,-0.75+0.125,-0.5+0.125,-0.25+0.125,0+0.125,0.25+0.125,0.5+0.125,0.75+0.125,1+0.125,1.25+0.125,1.5+0.125,1.75+0.125} { \pgfmathparse{0.020*rand+0.05} \let \z \pgfmathresult;
 \filldraw[draw=black,fill=gray!60] (\col,\hi) circle (\z);}}

\fill[fill=white, opacity=1] (0,0) .. controls (1/2,-1) and (1,1) ..  (2,1) .. controls (3.3,1.3) and (-1/2,3) .. (-1,1) -- (-1,2) -- (3.1,2) -- (3.1,-1) -- (-3,-1) -- (-0.1,-1) -- cycle;
\fill[fill=white, opacity=1] (-1,1) .. controls (-1.1,1/2) and (-0.3,0.3) .. (0,0) -- (0.1,-1) -- (-1.6,-1) -- (-1.6,2) -- (-0.9,2) -- cycle;
\draw (0,0) .. controls (1/2,-1) and (1,1) ..  (2,1) .. controls (3.3,1.3) and (-1/2,3) .. (-1,1) .. controls (-1.1,1/2) and (-0.3,0.3) .. (0,0);
\fill[fill=gray!60, opacity=1] (-1.5,-0.5) circle (0.3);
\draw (-1.5,-0.5) circle (0.3);
\draw (-2,-1) rectangle (-1,0);
\draw[dashed] (-1,0) -- (0,0.75);
\draw[dashed] (-1,-1) -- (0,0.50) node[midway,right]{$\times \e_n$};
\draw[dashed] (-1.5,-0.8) -- (-2.5,-0.8);
\draw[dashed] (-1.5,-0.8+0.6) -- (-2.5,-0.8+0.6);
\draw[<->] (-2.5,-0.8+0.6) -- (-2.5,-0.8) node[midway,left]{$2 \, r_{n,k}$};
\draw (2,-0.5) -- (2.15,-0.5) -- (2.15,-0.35) -- (2,-0.35) -- cycle;
\fill[fill=gray!60, opacity=1] (2,-0.5+0.4) -- (2.15,-0.5+0.4) -- (2.15,-0.35+0.4) -- (2,-0.35+0.4) -- cycle;
\draw (2,-0.5+0.4) -- (2.15,-0.5+0.4) -- (2.15,-0.35+0.4) -- (2,-0.35+0.4) -- cycle;
\draw (2,-0.5) -- (2.15,-0.5) -- (2.15,-0.35) -- (2,-0.35) -- cycle;
\draw (2,-0.5) -- (2.15,-0.5) -- (2.15,-0.35) -- (2,-0.35) -- cycle;
\draw (2.15,-0.39+0.4) -- (2.15,-0.39+0.4) node[midway,right]{$\w \Omega_n$};
\draw (2.15,-0.39) -- (2.15,-0.39) node[midway,right]{$\w \Omega \setminus \w \Omega_n$};
\end{tikzpicture}
\caption{The cross section of the non periodic microstructure}
\label{fignonpercirc}
\end{figure}
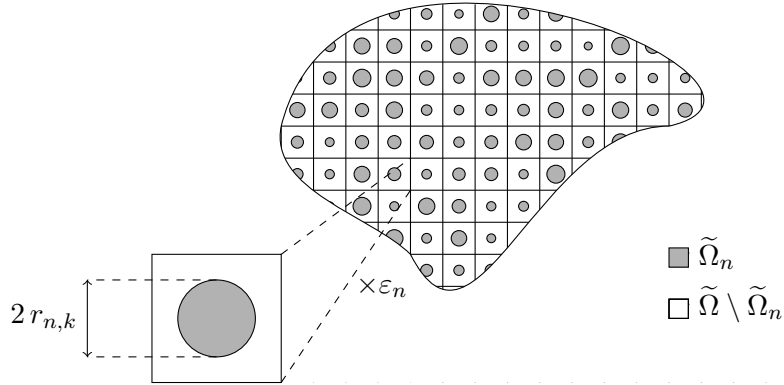

Let~${\rho}$ be a continuous function on~${\overline \Omega}$ depending only on the variable~${x' = (x_1,x_2)}$ satisfying 
\begin{equation} \label{thetapos}
\exists \, c_1, c_2 >0, \quad c_1 \leq \rho(x') \leq c_2, \quad \forall \, x = (x',x_3) \in \overline \Omega,
\end{equation}
and let~${r_n}$ be a sequence of positive numbers converging to~${0}$, as~${n}$ goes to infinity. We assume, without loss of generality, that 
\begin{equation} \label{condrho1}
\fint_\Omega \rho \dx = 1. 
\end{equation} 
Define, for any~${k \in \Z^3}$, the sequence~${(r_{n,k})_{n \in \N^*}}$ by 
\begin{equation*} 
{r_{n,k}:= r_n \, \sqrt{\rho(\e_n k)}.}
\end{equation*} 
We consider the case where~${\Omega_n}$ is the set of circular fibres ~${{\omega_{n,k} = \{ y \in Y \ | \ y_1^2 + y_2^2 \leq r_{n,k}^2\}}}$~(see Figure~\ref{fignonpercirc}) 
\begin{equation} \label{omegacirc}
\Omega_n = \Omega \cap \bigcup_{k \in \Z^3} \e_n (\omega_ {n,k} + k).
\end{equation}
Note that the fibres~${\omega_{n,k}}$ do not have the same radius.
\noindent We have the following result:
\begin{Prop} \label{thfib}
Let~${\Omega_n}$ be the sequence of subsets of~${\Omega}$ defined by~\eqref{omegacirc} and~${\sigma_n(h)}$ be the associated conductivity in the problem~\eqref{Dir}. Assume that
\begin{equation} \label{condcirc}
\e_n^2 |\ln r_n| \xrightarrow[n \to \infty]{} 0.
\end{equation}

\noindent Then, there exist a matrix-valued function~${\sigma_*(h)}$ and a subsequence of~${n}$, still denoted by~${n}$, such that the solution~${u_n}$ of the problem~\eqref{Dir} converges weakly in~${H_0^1(\Omega)}$ to the solution~${u}$ of the conductivity problem
\begin{equation*}
\left\{ \begin{array}{r c l l}
- \divg \big(\sigma_*(h) \nabla u \big) & = & f &\text{in } \Omega,\\
u & = & 0 & \text{on } \partial \Omega,
\end{array} \right.
\end{equation*} 
where~${\sigma_*(h)}$ is given, for any~${x = (x',x_3) \in \Omega}$, by 
\begin{equation} \label{eqperex1}
\sigma_*(h)(x) = \alpha_1 I_3 + \displaystyle \rho(x') \left(\frac{\alpha_2^3 + \alpha_2 \beta_2^2 |h|^2}{\alpha_2^2 + \beta_2^2 h_3^2} \right) e_3 \otimes e_3 + \beta_1 \mathscr{E}(h).
\end{equation}
\end{Prop}

\begin{Rem} \label{bidcirc} \upshape 
We can easily check that the homogenized conductivity~${\w \sigma_*(h)}$ of the two dimensional microstructure of the Figure~\ref{fignonpercirc} is given by 
\begin{equation*}  
\w \sigma_*(h) = \alpha_1 I_2 + \beta_1 h_3 J = \w \sigma_1(h).
\end{equation*} 
This leads to the simple form~\eqref{eqperex1} of~${\sigma_*(h)}$. 
\end{Rem}


\noindent \textbf{Proof of Proposition~\ref{thfib}} In order to apply Theorem~\ref{thgen}, we need to check that the conditions~\mbox{\eqref{WPW}-\eqref{ii}} hold true. On the one hand, the Poincaré-Wirtinger inequality combined with~\eqref{convphas20} imply the existence, for any~${k \in \Z^3}$, of a sequence of positive constants~${c_{n,k}}$ such that
\begin{equation} \label{WPWdem}
\forall \, v \in H^1(Q_n^k), \quad \int_{Q_n^k} \alpha_n \left| v - \fint_{Q_n^k} v\right|^2 \dx \leq c_{n,k} \int_{Q_n^k} \alpha_n |\nabla v|^2 \dx,
\end{equation} 
where~${Q_n^k = \e_n (Y + k)}$. 
Using estimates derived in~\cite{BrARMA}, one can show that the best constant in the weighted Poincaré-Wirtinger inequality~\eqref{WPWdem} satisfies
\begin{equation*}  
\forall \, k \in \Z^3, \quad \forall \, n \in \N^*, \quad 0< c_{n,k} \leq c \, \e_n^2 \left| \ln\big(r_n \, \sqrt{\rho(\e_n k)} \big)\right|,
\end{equation*} 
for some positive constant~${c}$. Therefore, by~\eqref{thetapos} and~\eqref{condcirc}, we have
\begin{equation*}
0 < c_{n,k} \leq c \, \e_n^2 |\ln r_n| + o(1) \xrightarrow[n \to \infty]{} 0, 
\end{equation*} 
uniformly with respect to~${ k \in \Z^3}$. Conditions~\eqref{WPW} and~\eqref{condWPW} of Theorem~\ref{thgen} are satisfied. On the other hand, by the definition of~${\Omega_n}$ and~\eqref{condrho1}, we have the following estimate for the volume fraction
\begin{equation*}
\theta_n = \frac{|\Omega_n|}{|\Omega|} \underset{n \to \infty}{\sim} \frac{1}{|\Omega|} \Sum_{\e_n k \in \Omega} \e_n^2 \pi r_n^2 \rho(\e_n k) \underset{n \to \infty}{\sim} \pi r_n^2 \fint_\Omega \rho \dx = \pi r_n^2, 
\end{equation*} 
which, by~\eqref{thetapos}, implies that for any~${n \in \N^*}$ and~${k \in \Z^3}$, 
\begin{equation*} 
\fint_{Q_n^k} \alpha_n \dx = \alpha_1 \big(1 - \pi r_n^2 \rho(\e_n k)\big) + \alpha_{2,n} \pi r_n^2 \rho(\e_n k) \leq c + c \, \theta_n \alpha_{2,n} \leq c.
\end{equation*} 
Then, condition~\eqref{ii} of Theorem~\ref{thgen} is satisfied. Theorem~\ref{thgen} and Remark~\ref{bidcirc} ensure the existence of an effective conductivity~${\sigma_*(h)}$ which, after an easy computation, writes 
\begin{equation} \label{eqperex1'}
\sigma_*(h)(x) = \alpha_1 I_3 + \displaystyle \theta(x') \left(\frac{\alpha_2^3 + \alpha_2 \beta_2^2 |h|^2}{\alpha_2^2 + \beta_2^2 h_3^2} \right) e_3 \otimes e_3 + \beta_1 \mathscr{E}(h) \quad \forall \, x = (x',x_3) \in \w \Omega \times (0,1),
\end{equation} 
where~${\theta}$ is the weak limit of~${\theta_n^{-1} \2}$. The function~${\theta}$ in~\eqref{eqperex1'} coincides with~${\rho}$. Indeed, since~${\rho}$ is continuous, we obtain, for any~${\varphi \in \mathscr C_0(\Omega)}$ extended to~${\R^3}$ by setting~${\varphi \equiv 0}$ on~${\R^3 \setminus \Omega}$, 
\begin{equation*} 
\int_\Omega \theta_n^{-1} \2 \varphi \dx = \frac{1}{\pi r_n^2} \Sum_{k \in \Z^3} \int_{\omega_{n,k}} \varphi \dx + o(1) = \frac{1}{\pi r_n^2} \Sum_{\e_n k \in \Omega} \e_n^2 \pi r_n^2 \rho(\e_n k ) \varphi(\e_n k) + o(1) 
\end{equation*} 
which implies that 
\begin{equation*} 
\int_\Omega \theta_n^{-1} \2 \varphi \dx = \int_\Omega \rho \varphi \dx + o(1).
\end{equation*} 
Finally~${\theta_n^{-1}\2}$ converges weakly-${*}$ to~${\rho}$ in~${\mathcal M(\Omega)}$ and, then,~${\theta \equiv \rho}$. This concludes the proof of Proposition~\ref{thfib}. 
\qed  

\subsection{Thin squared grids}

In this section, we consider the case of a columnar composite the cross section of which is a highly conducting grid surrounded by another conducting medium~(see Figure~\ref{fignonpersqa}). 
Let~${t_n}$ be a positive sequence converging to~${0}$ as~${n}$ goes to infinity. Let~${\rho}$ be a continuous function on~${\overline \Omega}$, depending only on the variable~${x' = (x_1,x_2)}$ and satisfying 
\begin{equation} \label{taupos}
\exists \, c_1, c_2 >0, \quad c_1 \leq \rho(x') \leq c_2, \quad \forall \, x = (x',x_3) \in \overline \Omega.
\end{equation}
We assume, without loss of generality, that 
\begin{equation*}
\fint_\Omega \rho \dx = 1. 
\end{equation*}  
Define, for any~${k}$ in~${\Z^3}$, the sequence~${(t_{n,k})_{n \in \N^*}}$ by 
\begin{equation} \label{condtn}
t_{n,k}:=\rho(\e_n k) \, t_n.
\end{equation} 
Let~${\Omega_n}$ be the set of non periodically distributed squared fibres 
\begin{equation} \label{crossper}
\Omega_n = \Omega \cap \bigcup_{k \in \Z^3} \e_n (\omega_ {n,k} + k) \quad \text{where }  \omega_{n,k} := \left\{ y \in Y \ | \ \max(|y_1|,|y_2|)\geq \textstyle{ \frac 1 2} -t_{n,k} \right\}.
\end{equation} 
Note that the case~${\rho \equiv 1}$ leads to a periodic distribution of the squared fibres in~${\Omega}$.\begin{center}
\begin{figure}[H]
\centering
\leavevmode
\hspace{-1.5cm}\subfloat[Periodic case]{\label{figdim3gene}
\begin{tikzpicture}[scale=1.7]
\fill[fill=gray!60, opacity=1] (-1.6,-0.9) rectangle (3,2);
\clip (-2.1,-1.1) rectangle (3,2);
\foreach \col in {-0.5,-0.25,0,0.25,0.5,0.75,1,1.25,1.5,1.75}{
\draw (-1.5,\col) -- (3,\col);}
\foreach \col in {-1.5,-1.25,-1,-0.75,-0.5,-0.25,0,0.25,0.5,0.75,1,1.25,1.5,1.75,2,2.25,2.5,2.75,3}{
\draw (\col,-1) -- (\col,1.95);}

\foreach \col in {-1.5+0.125,-1.25+0.125,-1+0.125,-0.75+0.125,-0.5+0.125,-0.25+0.125,0+0.125,0.25+0.125,0.5+0.125,0.75+0.125,1+0.125,1.25+0.125,1.5+0.125,1.75+0.125,2+0.125}{\foreach \hi in {-1+0.125,-0.75+0.125,-0.5+0.125,-0.25+0.125,0+0.125,0.25+0.125,0.5+0.125,0.75+0.125,1+0.125,1.25+0.125,1.5+0.125,1.75+0.125} {

\filldraw[draw=black,fill=white] (\col-0.08,\hi-0.08) rectangle (\col+0.08,\hi+0.08);}}
\fill[fill=white, opacity=1] (0,0) .. controls (1/2,-1) and (1,1) ..  (2,1) .. controls (3.3,1.3) and (-1/2,3) .. (-1,1) -- (-1,2) -- (3.1,2) -- (3.1,-1) -- (-3,-1) -- (-0.1,-1) -- cycle;

\fill[fill=white, opacity=1] (-1,1) .. controls (-1.1,1/2) and (-0.3,0.3) .. (0,0) -- (0.1,-1) -- (-1.6,-1) -- (-1.6,2) -- (-0.9,2) -- cycle;
\draw (0,0) .. controls (1/2,-1) and (1,1) ..  (2,1) .. controls (3.3,1.3) and (-1/2,3) .. (-1,1) .. controls (-1.1,1/2) and (-0.3,0.3) .. (0,0);

\fill[fill=gray!60, opacity=1] (-2+3.2525,-1) rectangle (-1+3.2525,0);
\fill[fill=white, opacity=1] (-2+3.2525+0.15,-1+0.15) rectangle (-1+3.2525-0.15,-0.15);
\draw[dashed] (-1+3.2525,-0.15) -- (-1+3.2525,0.5);
\draw[dashed] (-1+3.2525-0.15,-0.15) -- (-1+3.2525-0.15,0.5);
\draw[<->] (-1+3.2525-0.15,0.5) -- (-1+3.2525,0.5) node[midway,above]{$t_n$};
\draw (-2+3.2525+0.15,-1+0.15) rectangle (-1+3.2525-0.15,-0.15);
\draw (-2+3.2525,-1) rectangle (-1+3.2525,0);
\draw[dashed] (-1+3.2525-1,0) -- (0+3.2525-3,0.75);
\draw[dashed] (-1+3.2525-1,-1) -- (0+3.2525-3,0.50) node[midway,below]{$\hspace{-0.4cm}\times \e_n$};
\end{tikzpicture}}
\hspace{-1.5cm} \subfloat[Non periodic case]{\label{figdim3cell}
\begin{tikzpicture}[scale=1.7]
\fill[fill=gray!60, opacity=1] (-1.6,-0.9) rectangle (3,2);
\clip (-2.1,-1.1) rectangle (3,2);
\foreach \col in {-0.5,-0.25,0,0.25,0.5,0.75,1,1.25,1.5,1.75}{
\draw (-1.5,\col) -- (3,\col);}
\foreach \col in {-1.5,-1.25,-1,-0.75,-0.5,-0.25,0,0.25,0.5,0.75,1,1.25,1.5,1.75,2,2.25,2.5,2.75,3}{
\draw (\col,-1) -- (\col,1.95);}

\foreach \col in {-1.5+0.125,-1.25+0.125,-1+0.125,-0.75+0.125,-0.5+0.125,-0.25+0.125,0+0.125,0.25+0.125,0.5+0.125,0.75+0.125,1+0.125,1.25+0.125,1.5+0.125,1.75+0.125,2+0.125}{\foreach \hi in {-1+0.125,-0.75+0.125,-0.5+0.125,-0.25+0.125,0+0.125,0.25+0.125,0.5+0.125,0.75+0.125,1+0.125,1.25+0.125,1.5+0.125,1.75+0.125} { \pgfmathparse{0.025*rand+0.06} \let \z \pgfmathresult;

\filldraw[draw=black,fill=white] (\col-\z,\hi-\z) rectangle (\col+\z,\hi+\z);}}
\fill[fill=white, opacity=1] (0,0) .. controls (1/2,-1) and (1,1) ..  (2,1) .. controls (3.3,1.3) and (-1/2,3) .. (-1,1) -- (-1,2) -- (3.1,2) -- (3.1,-1) -- (-3,-1) -- (-0.1,-1) -- cycle;

\fill[fill=white, opacity=1] (-1,1) .. controls (-1.1,1/2) and (-0.3,0.3) .. (0,0) -- (0.1,-1) -- (-1.6,-1) -- (-1.6,2) -- (-0.9,2) -- cycle;
\draw (0,0) .. controls (1/2,-1) and (1,1) ..  (2,1) .. controls (3.3,1.3) and (-1/2,3) .. (-1,1) .. controls (-1.1,1/2) and (-0.3,0.3) .. (0,0);
\fill[fill=gray!60, opacity=1] (2,-0.5+0.4) -- (2.15,-0.5+0.4) -- (2.15,-0.35+0.4) -- (2,-0.35+0.4) -- cycle;
\draw (2,-0.5+0.4) -- (2.15,-0.5+0.4) -- (2.15,-0.35+0.4) -- (2,-0.35+0.4) -- cycle;
\draw (2,-0.5) -- (2.15,-0.5) -- (2.15,-0.35) -- (2,-0.35) -- cycle;
\draw (2,-0.5) -- (2.15,-0.5) -- (2.15,-0.35) -- (2,-0.35) -- cycle;
\draw (2.15,-0.39+0.4) -- (2.15,-0.39+0.4) node[midway,right]{$\w \Omega_n$};
\draw (2.15,-0.39) -- (2.15,-0.39) node[midway,right]{$\w \Omega \setminus \w \Omega_n$};
\end{tikzpicture}}
\caption{The cross section of the structure}
\label{fignonpersqa}
\end{figure}
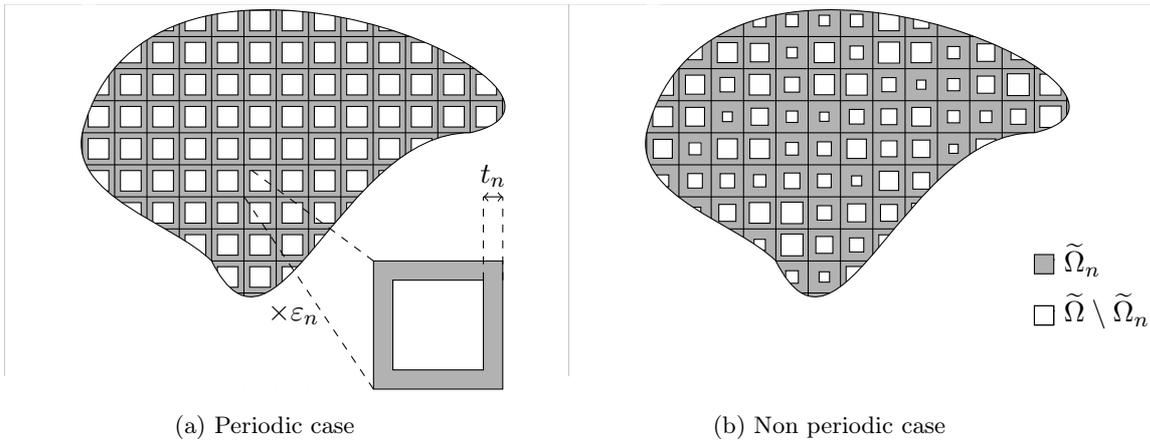
\end{center}

\noindent We have the following result:

\begin{Prop} \label{thcross}
Let~${\Omega_n}$ be the sequence of subsets of~${\Omega}$ defined by~\eqref{crossper} and~${\sigma_n(h)}$ be the associated conductivity in the problem~\eqref{Dir}. Assume that
\begin{equation} \label{condalphasim!}
4 \, t_n \alpha_{2,n} \xrightarrow[n \to \infty]{} \alpha_2 > 0 \ \ \textrm{and} \ \ 4 \, t_n \beta_{2,n} \xrightarrow[n \to \infty]{} \beta_2 \in \R.
\end{equation}

\noindent Then, there exist a matrix-valued function~${\sigma_*(h)}$ and a subsequence of~${n}$, still denoted by~${n}$, such that the solution~${u_n}$ of the problem~\eqref{Dir} converges weakly in~${H_0^1(\Omega)}$ to the solution~${u}$ of the conductivity problem
\begin{equation} \label{Dir*nonper2}
\left\{ \begin{array}{r c l l}
- \divg \big(\sigma_*(h) \nabla u \big) & = & f &\text{in } \Omega,\\
u & = & 0 & \text{on } \partial \Omega,
\end{array} \right. 
\end{equation}
where~$\sigma_*(h)$ is given by
\begin{equation} \label{propper'ex}
\sigma_*(h) := \begin{pmatrix}
\w \sigma_*(h) & p_* \\*[0.4em]
q_*^\T & \alpha_*
\end{pmatrix}
\end{equation} 
and, for any~${(x',x_3) \in \Omega}$, 
\begin{equation} \label{eqperex}
\left\{ \begin{array}{l}
\vspace{0.4cm}\displaystyle \w\sigma_*(h) = \left(\alpha_1 + \rho(x') \, \frac{\alpha_2^2 + \beta_2^2 h_3^2 }{2\alpha_2}\right)I_2 + \beta_1 h_3 J,\\
\vspace{0.2cm} \displaystyle p_* = - \left[\beta_1 + \rho(x') \, \frac{\beta_2}{2}\right] J \w h + \rho(x') \, \frac{\beta_2^2 \, h_3}{2 \alpha_2} \w h, \\
\vspace{0.2cm}q_* = \displaystyle \left[\beta_1 + \rho(x') \, \frac{\beta_2}{2}\right] J \w h + \rho(x') \, \frac{\beta_2^2\, h_3}{2 \alpha_2} \w h, \\
\alpha_* = \alpha_1 + \rho(x') \, \alpha_2 +   \rho(x') \, \displaystyle \frac{\beta_2^2}{2 \, \alpha_2} (h_1^2 + h_2^2).
\end{array}\right.
\end{equation}
\noindent In formula~\eqref{eqperex},~${\rho \equiv 1}$ corresponds to the periodic case.
\end{Prop}

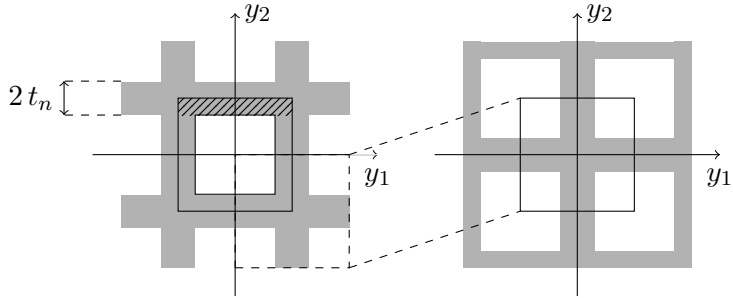
\begin{figure}[H]
\centering
\begin{tikzpicture}[scale=0.75,xmin=-18,xmax=13.5,ymin=-6,ymax=6]
\fill[fill=white, opacity=1] (-4,-3) rectangle (4,3);
\fill[fill=gray!60, opacity=1] (-2,-2) rectangle (2,2);
\fill[fill=white] (-1+0.3,-1+0.3) rectangle (1-0.3,1-0.3);
\fill[fill=white] (1+0.3,-1+0.3) rectangle (3-0.3,1-0.3);
\fill[fill=white] (-3+0.3,-1+0.3) rectangle (-1-0.3,1-0.3);
\fill[fill=white] (-1+0.3,-1+0.3+2) rectangle (1-0.3,1-0.3+2);
\fill[fill=white] (1+0.3,-1+0.3+2) rectangle (3-0.3,1-0.3+2);
\fill[fill=white] (-3+0.3,-1+0.3+2) rectangle (-1-0.3,1-0.3+2);
\fill[fill=white] (-1+0.3,-1+0.3-2) rectangle (1-0.3,1-0.3-2);
\fill[fill=white] (1+0.3,-1+0.3-2) rectangle (3-0.3,1-0.3-2);
\fill[fill=white] (-3+0.3,-1+0.3-2) rectangle (-1-0.3,1-0.3-2);
\draw (-1,-1) rectangle (1,1);
\draw (-1+0.3,-1+0.3) rectangle (1-0.3,1-0.3);
\draw[dashed] (0,-1-1) rectangle (2,1-1);
\draw[->] (-2.5,0) -- (2.5,0);
\draw[->] (0,-2.5) -- (0,2.5);
\draw (2.5,-0.4) node{$y_1$};
\draw (0.4,2.5) node{$y_2$};
\fill[white] (1-0.3+1.3,0) -- (4,0) ;
\draw[dashed] (-1+0.3-1.3,-1+0.3+2) -- (-3,-1+0.3+2) ;
\draw[dashed] (-1+0.3-1.3,-1+2-0.3) -- (-3,-1-0.3+2) ;
\draw[<->] (-3,-1-0.3+2) -- (-3,-1+0.3+2) node[midway,left]{$2 \, t_n$};
\fill[ pattern= north east lines] (-1,0.69) rectangle (1,1);

\fill[fill=gray!60, opacity=1] (-2+6,-0.3) rectangle (2+6,+0.3);
\fill[fill=gray!60, opacity=1] (-2+6,-0.3+2) rectangle (2+6,+2);
\fill[fill=gray!60, opacity=1] (-2+6,-2) rectangle (2+6,-2+0.3);
\fill[fill=gray!60, opacity=1] (-0.3+6,-2) rectangle (0.3+6,2);
\fill[fill=gray!60, opacity=1] (-0.3+6+2,-2) rectangle (6+2,2);
\fill[fill=gray!60, opacity=1] (6-2,-2) rectangle (0.3+6-2,2);
\draw (-1+6,-1) rectangle (1+6,1);
\draw[->] (-2.5+6,0) -- (2.5+6,0);
\draw[->] (0+6,-2.5) -- (0+6,2.5);
\draw (2.5+6,-0.4) node{$y_1$};
\draw (0.4+6,2.5) node{$y_2$};
\draw[dashed] (2,0) -- (2+3,1);
\draw[dashed] (2,0-2) -- (2+3,-1);
\end{tikzpicture}
\caption{The period cell of the cross section of the microstructure}
\label{figcrossperproof}
\end{figure}

\noindent \textbf{Proof of Proposition~\ref{thcross}.} Let us first consider the periodic case. In order to apply Proposition~\ref{propper}, we need to check that~\eqref{WPWper} and~\eqref{condWPWper} are satisfied. To this end, consider~${V \in \mathscr C^1(\overline{Y})}$ such that~${\langle V \rangle = 0}$. Define, for any~${n \in \N^*}$, the subsets~${K_{n}^i}$,~${i=1,2,3,4}$, of~${\overline{Y}}$ by 
\begin{equation*} 
\begin{array}{r l} K_{n}^1 := \left[ \textstyle{-\frac{1}{2}, \frac{1}{2}} \right] \times \left[ \textstyle{\frac{1}{2} - t_n, \frac{1}{2}} \right] \times \left[ \textstyle{-\frac{1}{2}, \frac{1}{2}} \right], & K_n^2 := \left[ \textstyle{-\frac{1}{2}, \frac{1}{2}} \right] \times \left[ \textstyle{-\frac{1}{2} + t_n, - \frac{1}{2}} \right] \times \left[ \textstyle{-\frac{1}{2}, \frac{1}{2}} \right], \\*[0.7em]
K_{n}^3 := \left[ \textstyle{\frac{1}{2} - t_n, \frac{1}{2}} \right] \times \left[ \textstyle{-\frac{1}{2}, \frac{1}{2}} \right] \times \left[ \textstyle{-\frac{1}{2}, \frac{1}{2}} \right], & K_n^4 := \left[ \textstyle{-\frac{1}{2} + t_n, - \frac{1}{2}} \right] \times \left[ \textstyle{-\frac{1}{2}, \frac{1}{2}} \right] \times \left[ \textstyle{-\frac{1}{2}, \frac{1}{2}} \right]. \end{array}
\end{equation*} 
For instance, the projection of~${K_n^1}$, in the~${(y_1,y_2)}$-plane, is the shaded zone in Figure~\ref{figcrossperproof}.

\noindent  By the definition~\eqref{defperan} of~${a_n}$, we have 
\begin{equation*} 
\begin{array}{r c l}
\displaystyle \int_Y a_n V^2 \dy & \leq & \displaystyle \int_{Y \setminus \overset{4}{\underset{i=1}{\cup}} K_n^i} a_n V^2 \dy + \displaystyle \Sum_{i=1}^4 \int_{K_{n}^i} a_n V^2 \dy \\
& \leq & \displaystyle \alpha_1 \int_{Y \setminus \overset{4}{\underset{i=1}{\cup}} K_n^i} V^2 \dy + \alpha_{2,n} \displaystyle \Sum_{i=1}^4 \int_{K_{n}^i} V^2 \dy \\
& \leq & \displaystyle \alpha_1 \int_{Y} V^2 \dy + \alpha_{2,n} \displaystyle \Sum_{i=1}^4 \int_{K_{n}^i} V^2 \dy.
\end{array}
\end{equation*}
Since~${\langle V \rangle=0}$, this inequality and the Poincaré-Wirtinger inequality in~${H^1(Y)}$, yield
\begin{equation} \label{ineqPW}
\displaystyle \int_Y a_n V^2 \dy \leq \displaystyle \alpha_1 \int_{Y} |\nabla V|^2 \dy + \alpha_{2,n} \displaystyle \Sum_{i=1}^4 \int_{K_{n}^i} V^2 \dy.
\end{equation}
We now estimate the second term of the right hand side of this inequality. On the one hand, since~${K_n^1}$ is convex, the Poincaré-Wirtinger constant in~${H^1(K_n^1)}$  is bounded from above by the diameter of~${K_n^1}$ divided by~${\pi}$~\cite{PaWe} and, therefore 
\begin{equation} \label{ineqPW'}
\begin{array}{r c l}
\displaystyle \int_{K_{n}^1} V^2 \dy & \leq & 2 \displaystyle \int_{K_{n}^1} \Big| V - \fint_{K_{n}^1} V \dy \Big|^2 + 2 \, |K_n^1| \ \Big|\fint_{K_{n}^1} V \dy \Big|^2 \\*[1.1em]
& \leq & \displaystyle c \left( \int_{ K_n^1} |\n V | ^2 \dy + |K_n^1| \ \Big|\fint_{K_{n}^1} V \dy \Big|^2 \right).
\end{array}
\end{equation} 
On the other hand, for any~${\textstyle{- \frac 1 2} \leq r, \, s, \, t \leq \textstyle{\frac 1 2}}$, we have 
\begin{equation} \label{ineqPW''}
\widehat V(s) - \widehat V(r) = \int_{r}^{s} \displaystyle \widehat V'(t) \ \text{d}t, \quad \text{where} \quad \widehat V(t) := \int_{-\frac{1}{2}}^{\frac{1}{2}} \int_{-\frac{1}{2}}^{\frac{1}{2}}  V(y_1,t,y_3) \ \text{d}y_1\text{d}y_3. 
\end{equation} 
Integrating the first equality in~\eqref{ineqPW''} with respect to~${s \in [ \textstyle{\frac 1 2} - t_n, \textstyle{\frac 1 2}]}$ and~${r \in [ \textstyle{- \frac 1 2}, \textstyle{\frac 1 2}]}$, we have 
\begin{equation*} 
\bigg|\fint_{K_{n}^1} V \dy - \fint_Y V \dy \bigg| \leq \fint_{\frac 1 2 -t_n}^{\frac 1 2} \fint_{- \frac 1 2}^{\frac 1 2}\int_{- \frac 1 2}^{\frac 1 2} |\widehat V'(t)| \ \text{d}t \leq \int_Y \left|\frac{\partial V}{\partial y_2} \right| \dy,
\end{equation*}
which, since~${\langle V \rangle = 0}$, implies that
\begin{equation} \label{ineqPW21}
\bigg|\fint_{K_{n}^1} V \dy \bigg|  \leq \int_Y |\nabla V| \dy \leq \Vert  \n V \Vert_{L^2(Y)^3}.
\end{equation} 
Then, combining~\eqref{ineqPW'} and~\eqref{ineqPW21} with the boundedness~\eqref{condalphasim!} of~${|K_n^1| \alpha_{2,n} =  t_n \alpha_{2,n}}$, we obtain that 
\begin{equation} \label{ineqPW'''} 
\alpha_{2,n} \int_{K_{n}^1} V^2 \dy \leq c \left(\alpha_{2,n} \int_{ K_n^1} |\n V | ^2 \dy + \Vert  \n V \Vert^2_{L^2(Y)^3} \right) \leq c \int_Y a_n |\n V|^2 \dy. 
\end{equation} 
Similarly to~\eqref{ineqPW'''}, we have, for~${i=2,3,4}$, 
\begin{equation} 
\alpha_{2,n} \int_{K_{n}^i} V^2 \dy \leq c \int_Y a_n |\n V|^2 \dy.\label{ineqPW'''1} 
\end{equation} 
Finally,~\eqref{ineqPW},~\eqref{ineqPW'''} and~\eqref{ineqPW'''1} imply 
\begin{equation} \label{per''''}
\int_Y a_n V^2 \dy \leq c \int_Y a_n | \n V|^2 \dy.
\end{equation} 
By a density argument,~\eqref{per''''} is satisfied for any~${V \in H^1(Y)}$ with~${\langle V \rangle = 0}$. Since~${\e_n}$ converges to~${0}$, the hypotheses~\eqref{WPWper} and~\eqref{condWPWper} of Proposition~\ref{propper} are satisfied. Then, there exists a homogenized matrix which is given in terms of the transversal effective conductivity~${\w \sigma_*}$ of the microstructure of Figure~\ref{figdim3gene}. It remains to determine
~${\w \sigma_*}$. Since one can choose the cross-like shape of the Figure~\ref{figcrossperproof} as the period cell of the transversal microstructure of the heterogeneous medium occupying~${\Omega}$, 
Proposition 3.2 of~\cite{BP1} ensures that
\begin{equation} \label{sigmat0} 
\w \sigma_* = \displaystyle \left(\alpha_1 + \frac{\alpha^2_2 + \beta^2_2 h_3^2}{2 \alpha_2}\right) I_2  + \beta_1 h_3 J,
\end{equation} 
and formula~\mbox{\eqref{propper'ex}-\eqref{eqperex}} of~${\sigma_*(h)}$ is a consequence of~\mbox{\eqref{propper'}-\eqref{eqper}} where~${\w \sigma_*}$ is given by~\eqref{sigmat0}. The periodic case is then proved.

\medskip

The existence of~${\sigma_*(h)}$ in the non periodic case is a consequence of Theorem~\ref{thgen}. Indeed, for any~${k \in \Z^3}$ and~${n\in\N^\ast}$, a rescaling of~\eqref{per''''} gives
\begin{equation*} 
 \forall \, v \in H^1(Q_n^k), \quad \int_{Q_n^k} \alpha_n \left| v - \fint_{Q_n^k} v\right|^2 \dx \leq c \, \e_n^2 \int_{Q_n^k} \alpha_n |\nabla v|^2 \dx,
\end{equation*} 
and, by~\eqref{condtn}, 
\begin{equation*} 
\fint_{Q_n^k} \alpha_n \dx = \alpha_1 \big(1 - 4 \, t_{n,k} (1-t_{n,k})\big) + 4 \, \alpha_{2,n} \, t_{n,k} (1-t_{n,k}) \leq c + c \, t_n \alpha_{2,n} \leq c.
\end{equation*}   
The assumptions of Theorem~\ref{thgen} are satisfied. Then, there exist a matrix-valued function~${\sigma_*(h)}$ and a subsequence of~${n}$, still denoted by~${n}$, such that the solution~${u_n}$ of the problem~\eqref{Dir} converges weakly in~${H_0^1(\Omega)}$ to the solution~${u}$ of the conductivity problem~\eqref{Dir*nonper2}. In view of the formulas~\eqref{propgen'} and~\eqref{eqgen}, the expression of~${\sigma_*(h)}$ becomes explicit as soon as~${\w \sigma_*}$ and~${\theta}$ are identified. 

\medskip

On the one hand, it is easy to check, similarly to the proof of Proposition~\ref{thfib}, that~${\rho}$ is the \mbox{weak-$*$} limit, in the sense of Radon measures, of the sequence~${\theta_n^{-1} \2}$. Then, the function~${\theta}$ in Theorem~\ref{thgen} turns out to be~${\rho}$. On the other hand, by Remark~\ref{remexpl'}, in order to compute~${\w \sigma_*}$, one has to determine~$\sigma_*^0(\alpha_1,\alpha_2)$, which is the~\mbox{${H(\mathcal M(\w \Omega)^2)}$-limit}, in the sense of Definition~\ref{HM}, of the conductivity~$\w \sigma_n(0)$, in the absence of a magnetic field, given by, for any~${x' \in \w \Omega}$,
\begin{equation*}
\w \sigma_n(0):=\left\{ \begin{array}{l l}
\alpha_1 I_2 & \text{in} \quad \w \Omega \setminus \w \Omega_n, \\
\alpha_{2,n} I_2 & \text{in} \quad \w \Omega_n.
\end{array}\right.
\end{equation*}
Due to the local nature~\cite{BrCa1} of the~\mbox{${H(\mathcal M(\w \Omega)^2)}$-convergence}, it is sufficient to compute~$\sigma_*^0(\alpha_1,\alpha_2)$ locally in~${\w \Omega}$. 
To that aim, consider~${x' \in \w \Omega}$ and~${\e > 0}$ small enough such that the closed disk~${\overline{D}(x',\e) \subset \w \Omega}$. Since~${\rho}$ is continuous and by~\eqref{taupos}, we have
\begin{equation} \label{eeee}
0 < c_{1,\e}(x') := \inf \limits_{z \in \overline{D}(x',\e)} \rho(z) \leq \rho(x') \leq c_{2,\e}(x') := \sup \limits_{z \in \overline{D}(x',\e)} \rho(z).
\end{equation} \vspace{-1.1cm}

\begin{center}
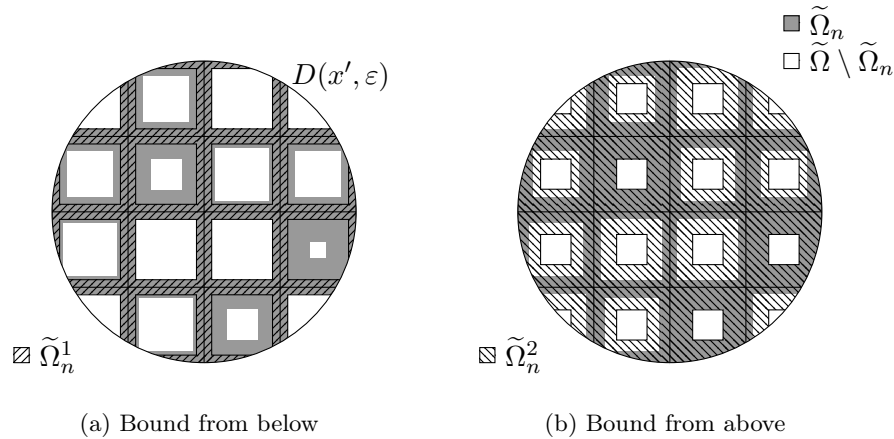
\begin{figure}[H]
\centering
\leavevmode
\subfloat[Bound from below]{
\begin{tikzpicture}[scale=1]
\clip (-3,-2.5) rectangle (3,3);
\fill[fill=gray!80, opacity=1] (-4,-4) rectangle (4,4);

\foreach \col in {-16,...,16}{
\draw (-4,\col) -- (4,\col);}
\foreach \col in {-16,...,16}{
\draw (\col,-4) -- (\col,4);}
\fill[white] (-1+0.1,-1+0.1) rectangle (0-0.1,0-0.1);
\fill[white] (-1+0.3,-1+0.3+1) rectangle (0-0.3,0-0.3+1);
\fill[white] (-1+0.2,-1+0.2+2) rectangle (0-0.2,0-0.2+2);
\fill[white] (-1+0.15,-1+0.15-1) rectangle (0-0.15,0-0.15-1);
\fill[white] (-1+0.1+1,-1+0.1) rectangle (0-0.1+1,0-0.1);
\fill[white] (-1+0.15+1,-1+0.15+1) rectangle (0-0.15+1,0-0.15+1);
\fill[white] (-1+0.1+1,-1+0.1+2) rectangle (0-0.1+1,0-0.1+2);
\fill[white] (-1+0.3+1,-1+0.3-1) rectangle (0-0.3+1,0-0.3-1);
\fill[white] (-1+0.4+2,-1+0.4) rectangle (0-0.4+2,0-0.4);
\fill[white] (-1+0.2+2,-1+0.2+1) rectangle (0-0.2+2,0-0.2+1);
\fill[white] (-1+0.1+2,-1+0.1+2) rectangle (0-0.1+2,0-0.1+2);
\fill[white] (-1+0.1-1,-1+0.1-1) rectangle (0-0.1-1,0-0.1-1);
\fill[white] (-1+0.15-1,-1+0.15) rectangle (0-0.15-1,0-0.15);
\fill[white] (-1+0.2-1,-1+0.2+1) rectangle (0-0.2-1,0-0.2+1);
\fill[white] (-1+0.1-1,-1+0.1+2) rectangle (0-0.1-1,0-0.1+2);
\fill[white] (-1+0.1+2,-1+0.1-1) rectangle (0-0.1+2,0-0.1-1);

\foreach \col in {-16,...,16}{
\fill[pattern= north east lines] (-4,\col-0.1) rectangle (4,\col+0.1);}
\foreach \col in {-16,...,16}{
\fill[pattern= north east lines] (\col-0.1,-4) rectangle (\col+0.1,4);}

\foreach \col in {-16,...,16}{
\foreach \h in {-16,...,16}{
\draw (-1+\col+0.1,-1+\h+0.1) rectangle (\col-0.1,\h-0.1);}}

\fill[white] (-0.1,-4) -- (0,-2) arc (270:450:2) -- (-0.1,4) -- (4,4) -- (4,-4) -- cycle;
\fill[white] (0.1,4) -- (0,2) arc (90:270:2) -- (0.1,-4) -- (-4,-4) -- (-4,4) -- cycle;
\draw (0,0) circle (2);
\draw (1.8,1.8) -- (1.8,1.8) node{${D(x',\e)}$};
\fill[pattern = north east lines] (-2.5,-2) -- (-2.3,-2) -- (-2.3,-1.8) -- (-2.5,-1.8) -- cycle;
\draw (-2.5,-2) -- (-2.3,-2) -- (-2.3,-1.8) -- (-2.5,-1.8) -- cycle;
\draw (-2.3,-1.9) -- (-2.3,-1.9) node[midway,right]{$\w \Omega_n^1$};
\end{tikzpicture}}
 \subfloat[Bound from above]{
\begin{tikzpicture}[scale=1]
\clip (-3,-2.5) rectangle (3,3);
\fill[fill=gray!80, opacity=1] (-4,-4) rectangle (4,4);

\foreach \col in {-16,...,16}{
\draw (-4,\col) -- (4,\col);}
\foreach \col in {-16,...,16}{
\draw (\col,-4) -- (\col,4);}
\fill[white] (-1+0.1,-1+0.1) rectangle (0-0.1,0-0.1);
\fill[white] (-1+0.3,-1+0.3+1) rectangle (0-0.3,0-0.3+1);
\fill[white] (-1+0.2,-1+0.2+2) rectangle (0-0.2,0-0.2+2);
\fill[white] (-1+0.15,-1+0.15-1) rectangle (0-0.15,0-0.15-1);
\fill[white] (-1+0.1+1,-1+0.1) rectangle (0-0.1+1,0-0.1);
\fill[white] (-1+0.15+1,-1+0.15+1) rectangle (0-0.15+1,0-0.15+1);
\fill[white] (-1+0.1+1,-1+0.1+2) rectangle (0-0.1+1,0-0.1+2);
\fill[white] (-1+0.3+1,-1+0.3-1) rectangle (0-0.3+1,0-0.3-1);
\fill[white] (-1+0.3+2,-1+0.3) rectangle (0-0.3+2,0-0.3);
\fill[white] (-1+0.2+2,-1+0.2+1) rectangle (0-0.2+2,0-0.2+1);
\fill[white] (-1+0.1+2,-1+0.1+2) rectangle (0-0.1+2,0-0.1+2);
\fill[white] (-1+0.1-1,-1+0.1-1) rectangle (0-0.1-1,0-0.1-1);
\fill[white] (-1+0.15-1,-1+0.15) rectangle (0-0.15-1,0-0.15);
\fill[white] (-1+0.2-1,-1+0.2+1) rectangle (0-0.2-1,0-0.2+1);
\fill[white] (-1+0.1-1,-1+0.1+2) rectangle (0-0.1-1,0-0.1+2);
\fill[white] (-1+0.1+2,-1+0.1-1) rectangle (0-0.1+2,0-0.1-1);

\foreach \col in {-16,...,16}{
\fill[pattern= north west lines] (-4,\col-0.3) rectangle (4,\col+0.3);}
\foreach \col in {-16,...,16}{
\fill[pattern= north west lines] (\col-0.3,-4) rectangle (\col+0.3,4);}

\foreach \col in {-16,...,16}{
\foreach \h in {-16,...,16}{
\draw (-1+\col+0.3,-1+\h+0.3) rectangle (\col-0.3,\h-0.3);}}


\fill[white] (-0.1,-4) -- (0,-2) arc (270:450:2) -- (-0.1,4) -- (4,4) -- (4,-4) -- cycle;
\fill[white] (0.1,4) -- (0,2) arc (90:270:2) -- (0.1,-4) -- (-4,-4) -- (-4,4) -- cycle;
\draw (0,0) circle (2);
\fill[fill=gray!80, opacity=1] (2-0.5,-0.5+0.4+2.5) -- (2.2-0.5,-0.5+0.4+2.5) -- (2.2-0.5,-0.3+0.4+2.5) -- (2-0.5,-0.3+0.4+2.5) -- cycle;
\draw (2-0.5,-0.5+0.4+2.5) -- (2.2-0.5,-0.5+0.4+2.5) -- (2.2-0.5,-0.3+0.4+2.5) -- (2-0.5,-0.3+0.4+2.5) -- cycle;
\draw (2.2-0.5,-0.39+0.4+2.5) -- (2.2-0.5,-0.39+0.4+2.5) node[midway,right]{$\w \Omega_n$};
\draw (2-0.5,-0.5+0.4+2.5-0.5) -- (2.2-0.5,-0.5+0.4+2.5-0.5) -- (2.2-0.5,-0.3+0.4+2.5-0.5) -- (2-0.5,-0.3+0.4+2.5-0.5) -- cycle;
\draw (2.2-0.5,-0.39+0.4+2) -- (2.2-0.5,-0.39+0.4+2) node[midway,right]{$\w \Omega \setminus \w \Omega_n$};
\fill[pattern = north west lines] (-2.5,-2) -- (-2.3,-2) -- (-2.3,-1.8) -- (-2.5,-1.8) -- cycle;
\draw (-2.5,-2) -- (-2.3,-2) -- (-2.3,-1.8) -- (-2.5,-1.8) -- cycle;
\draw (-2.3,-1.9) -- (-2.3,-1.9) node[midway,right]{$\w \Omega_n^2$};
\end{tikzpicture}}
\caption{Bounds from below and above of~${\w \sigma_n(0)}$}
\label{derniere}
\end{figure}
\end{center} \vspace{-1cm}

\noindent For~${i=1,2}$, let~${\w \Omega_n^{i}}$ be the subset of~${\overline D(x',\e)}$ defined by (see Figure~\ref{derniere}) 
\begin{equation} \label{crossperhhc}
\w \Omega_n^{i,\e} = \overline D(x',\e) \cap \bigcup_{k \in \Z^2} \e_n \left(k +  \left\{ y \in (\textstyle{-\frac 1 2}, \textstyle{\frac 1 2})^2 \ | \ \max(|y_1|,|y_2|)\geq \textstyle{\frac 1 2}-c_{i,\e}(x') \, t_n \right\}\right),
\end{equation}
and let~${\w \sigma_n^i}$ be the periodic conductivity defined on~${\overline D(x',\e)}$ by
\begin{equation} \label{condencad}
\w \sigma_n^{i,\e}:=\left\{ \begin{array}{l l}
\alpha_1 I_2 & \text{in} \quad \overline D(x',\e) \setminus \w \Omega_n^{i,\e}, \\
\alpha_{2,n} I_2 & \text{in} \quad \w \Omega_n^{i,\e}.
\end{array}\right.
\end{equation}
By the definitions~\eqref{eeee} and~\eqref{condencad}, we have for any~${z \in \overline D(x',\e)}$, the inequalities
\begin{equation} \label{ffff}
\w \sigma_n^{1,\e}(z) \leq \w \sigma_n(0)(z) \leq \w \sigma_n^{2,\e}(z).
\end{equation} 
For the rest of the proof, we need the following result which is a consequence of the two-dimensional div-curl lemma, in a high contrast context, of~\cite{BrCa1}:

\begin{Lem} \label{lemineq} Let~${D}$ be a bounded domain of~${\R^2}$  and, for~${i = 1,2}$, consider an equi-coercive sequence of symmetric matrix-valued functions~${A_n^i \in L^\infty(D)^{2 \times 2}}$ bounded in~${L^1(D)^{2 \times 2}}$ which~\mbox{${H(\mathcal M(D)^2)}$}-con\-verges to~${A_*^i}$ in the sense of Definition~\ref{HM}. We assume that
\begin{equation} \label{ineqlem22}
\forall \, n \in \N^*, \quad A_n^1 \leq A_n^2 \quad \text{a.e. in } D.
\end{equation} 
Then, we have the inequality 
\begin{equation*} 
 A_*^1 \leq A_*^2 \quad \text{a.e. in } D.
\end{equation*}
\end{Lem}
%

\noindent \textbf{Proof of Lemma~\ref{lemineq}} Let~${\lambda \in \R^2}$. Consider, for~${i=1,2}$, the solution~${v_n^{\lambda,i}}$ of 
\begin{equation*}
\left\{\begin{array}{r c l l}
\divg \big(A_n^i \w \nabla v_n^{\lambda,i} \big) & = & \mathrm{div} \big(A_*^i\lambda\big)\vspace{0.2cm} &\text{in } D,\\
v_n^{\lambda,i} & = & \lambda \cdot x & \text{on } \partial D.\\
\end{array} \right.
\end{equation*} 
By Definition~\ref{HM}, we have the convergences, for~${i=1,2}$,  
\begin{equation*}
\left\{ \begin{array}{ r l l}
v_n^{\lambda,i} & \harpoon \lambda \cdot x \vspace{0.2cm} & \text{weakly in } H^1(D), \\
A_n^i \w \nabla v_n^{\lambda,i} & \harpoon A_*^i\lambda &\text{weakly-${*}$ in }\mathcal M (D)^2.\\
\end{array} \right.
\end{equation*} 
On the one hand, by~\eqref{ineqlem22}, we have the inequality, almost everywhere in~${D}$
\begin{equation} \label{lem22}
2 \, A_n^1 \w \n v_n^{\lambda,1} \cdot \w \n v_n^{\lambda,2} - A_n^1 \w \n v_n^{\lambda,1} \cdot \w \n v_n^{\lambda,1} \leq A_n^1 \w \n v_n^{\lambda,2} \cdot \w \n v_n^{\lambda,2} \leq A_n^2 \w \n v_n^{\lambda,2} \cdot \w \n v_n^{\lambda,2}.
\end{equation}
On the other hand, applying, for~${i,j=1,2}$, the two-dimensional div-curl lemma of~\cite{BrCa1}~(Theorem~2.1) to~${\xi_n:= A_n^i \w \n v_n^{\lambda,i}}$ and~${v_n:=v_n^{\lambda,j}}$, we have the convergences, in the sense of distributions,
\begin{equation}  \label{lem11}
\forall \, i,j=1,2, \quad \xi_n \cdot \w \n v_n = A_n^i \w \n v_n^{\lambda,i} \cdot \w \n v_n^{\lambda,j} \harpoon A_*^i \lambda \cdot \lambda \quad \text{in} \quad \mathscr D'(D).
\end{equation}
Finally, combining~\eqref{lem22} and~\eqref{lem11}, we obtain
\begin{equation*}
2 \, A_*^1\lambda \cdot \lambda - A_*^1\lambda \cdot \lambda \leq A_*^2\lambda \cdot \lambda \quad \text{in} \quad \mathscr D'(D),  
\end{equation*}
which concludes the proof of Lemma~\ref{lemineq}.
\qed

\medskip

\noindent Since, for~${i=1,2}$,~${\w \sigma_n^{i,\e}}$ is an equi-coercive sequence of periodic matrix-valued functions bounded in~${L^1(D(x',\e))}$,~${\w \sigma_n^{i,\e}}$~\mbox{${H(\mathcal M(D)^2)}$-converges} to a constant matrix~${\w \sigma_*^{i,\e}}$. Then, applying Lemma~\ref{lemineq} with~${D = D(x',\e)}$ and~\eqref{ffff}, we have
\begin{equation} \label{eeeee}
\w \sigma_*^{1,\e} \leq \sigma_*^0(\alpha_1,\alpha_2)(z) \leq \w \sigma_*^{2,\e} \quad \text{ a.e.} \; z \in \overline D(x',\e).
\end{equation}
Moreover, due to the definition~\eqref{crossperhhc} of~${\w \Omega_n^{i,\e}}$ and the convergence~\eqref{condalphasim!}, we have
\begin{equation*}
|\w \Omega_n^{i,\e}|\, |\w \Omega|^{-1} \alpha_{2,n} = 4\, t_n \alpha_{2,n} \, c_{i,\e}(x') + o(1) \xrightarrow[n \to \infty]{} \alpha_{2} \, c_{i,\e}(x') >0.
\end{equation*}
Then, substituting~${\alpha_{2} \, c_{i,\e}(x')}$ for~${\alpha_2}$ in~\eqref{sigmat0} in the absence of a magnetic field~(\textit{i.e.},~${h_3 = 0}$), we obtain, for~$i = 1,2$,
\begin{equation} \label{eeee'}
\w \sigma_*^{i,\e} = \left(\alpha_1 + c_{i,\e}(x') \, \frac{\alpha_2}{2}  \right) I_2.
\end{equation}
By~\eqref{eeee} and~\eqref{eeee'}, taking the limit, as~${\e}$ goes to~${0}$, in the inequalities~\eqref{eeeee}, we obtain, for any Lebesgue point~${x'}$ of~${\sigma_*^0(\alpha_1,\alpha_2)}$ in~${ \w \Omega}$,
\begin{equation*} 
\sigma_*^0(\alpha_1,\alpha_2) = \left(\alpha_1 + \rho(x') \, \frac{\alpha_2}{2}\right) I_2.
\end{equation*}
Therefore, by Remark~\ref{remexpl'}, we have
\begin{equation} \label{tttt'}
\w \sigma_* = \left(\alpha_1 + \rho(x') \, \frac{\alpha_2^2 + \beta_2^2 h_3^2}{2 \, \alpha_2} \right) I_2 + \beta_1 h_3 J.
\end{equation}
Finally, we apply the formula~\mbox{\eqref{propgen'}-\eqref{eqgen}} for~${\sigma_*(h)}$ in Theorem~\ref{thgen}, with~${\w \sigma_*}$ given by~\eqref{tttt'}, to obtain~\mbox{\eqref{propper'ex}-\eqref{eqperex}}. This concludes the proof of Proposition~\ref{thcross}.\qed

\bigskip

\noindent \textit{Acknowledgements:} The authors wish to thank M. Briane for his valuable comments and suggestions. This work was partially supported by the ANR HJnet~(ANR-12-BS01-0008-01).

\markboth{Chapitre~\ref*{chap2} - Homogenization of high-contrast and non symmetric conductivities for non periodic columnar structures}{ \textit{References}}
\renewcommand\bibname{References}
\markboth{Chapitre~\ref*{chap2} - Homogenization of high-contrast and non symmetric conductivities for non periodic columnar structures}{ \textit{References}}
\bibliographystyle{plain}
\bibliography{biblio}
\end{document}